\documentclass[11pt]{article}
\usepackage{amsmath}
\usepackage{amsfonts}
 \newcounter{mnotecount}[section]

 \renewcommand{\themnotecount}{\thesection.\arabic{mnotecount}}

 \newcommand{\mnote}[1]
 {\protect{\stepcounter{mnotecount}}$^{\mbox{\footnotesize
 $
 \bullet$\themnotecount}}$ \marginpar{
 \raggedright\tiny\em
 $\!\!\!\!\!\!\,\bullet$\themnotecount: #1} }

\def\p{\partial}
\def\be{\begin{equation}}
\def\ee{\end{equation}}
\def\bea{\begin{eqnarray}}
\def\eea{\end{eqnarray}}
\def\bean{\begin{eqnarray*}}
\def\eean{\end{eqnarray*}}

\def\p{\partial}

\begin{document}

\title{Rod Structures and Patching Matrices: a review}
\author{Paul Tod\footnote{email: tod@maths.ox.ac.uk }\\Mathematical Institute,\\Oxford University}

\maketitle
\begin{abstract}
I review the twistor theory construction of stationary and axisymmetric, Lorentzian signature solutions of the Einstein vacuum equations and the related toric Ricci-flat metrics of Riemannian signature, \cite{W,MW,F,FW}. The construction arises from combining the Ward construction \cite{W2} of anti-self-dual Yang-Mills fields as holomorphic vector bundles on twistor space, with the observation of Witten \cite{LW} that the Einstein equations for these metrics include the anti-self-dual Yang-Mills equations. The principal datum for a solution is the holomorphic patching matrix $P$ for a rank-2 holomorphic vector bundle on a reduced twistor space, and $P$ is typically much simpler than the corresponding metric to write down.

I give a catalogue of examples, building on earlier collections \cite{F,AG}, and consider two inverse problems: how does the rod structure of such a metric, together with its asymptotics, determine $P$? And hoes does $P$ fix the metric?

\end{abstract}

\section{Introduction}
Given the fairly recent discovery of the Chen-Teo metric \cite{CT1,CT} and the recent classification of Hermitian toric ALF gravitational instantons \cite{BG}, this seems a good moment to offer a review of the twistor construction of four-dimensional Ricci-flat metrics with two commuting symmetries and to collect together the known examples. A word on terminology: the method can be adapted to Lorentzian or Riemannian signature; in the Lorentzian case the symmetry is usually described as `stationarity and axisymmetry' (though other two-dimensional symmetry groups may be of interest and can be handled), because the time-symmetry isn't made periodic, while in the Riemannian case it's always `toric', i.e. with two periodic commuting symmetries.

Historically, the subject begins with Lou Witten's paper \cite{LW} which noted that the equation governing static, axisymmetric self-dual\footnote{It's not important whether we specify self-dual or anti-self-dual here, and the former is shorter to write!} $SU(2)$ Yang-Mills fields in Euclidean 4-space was the same as the main equation for stationary, axisymmetric solutions of the Einstein vacuum equations in Lorentzian 4-space. This was picked up by Ward \cite{W} who observed that his construction \cite{W2} of solutions to self-dual Yang-Mills theory in terms of holomorphic vector bundles on (regions of) flat twistor space $\mathbb{CP}^3$, could therefore be applied to generate these Einstein vacua, giving solutions in terms of holomorphic vector bundles over a twistor space reduced by the symmetry. 

The next step was taken by Mason and Woodhouse \cite{MW}. They gave a detailed account of the reduced twistor space, which is essentially the twistor space $\mathbb{CP}^3$ or a large region in it, quotiented by two holomorphic vector fields derived from the two space-time symmetries. The action is not free and proper and the reduced twistor space $R$ turns out to be a one-dimensional non-Hausdorff complex manifold, essentially two copies of the Riemann sphere $\mathbb{CP}^1$ with identifications. It's still the case that holomorphic vector bundles $E$ on $R$ can be defined; this is done most conveniently with respect to a fixed standard four-set cover of $R$, two sets for each 
$\mathbb{CP}^1$, and then just one of the transition matrices has the information of the space-time in it -- call this the \emph{patching matrix} $P$ --  and the information of the metric can be extracted by splitting $P$ according to Ward's prescription \cite{W}. Different reality conditions lead to real Riemannian or real Lorentzian solutions.

Fletcher and Woodhouse \cite{F,FW} gave more examples and considered the inverse problem: without looking at the metric, what data from the space-time determine $P$ and therefore the metric? In this way, the Lorentzian Kerr metric is determined by its \emph{rod structure} (defined below) and its asymptotic conserved quantities. Clearly this gives another route to black hole uniqueness, which was their aim, and emphasises the status of $P$ as the fundamental datum for the metric. Finally, I should mention the extension of the method to five-dimensional black holes with three commuting symmetries due to Metzner and coworkers \cite{NM1,NM2} and the catalogue of examples compiled by Gray \cite{AG}.

The contents of this article are as follows: in Section 2, I briefly review the theory of stationary, axisymmetric space-times and give a sketch of the twistor theory described above; the notion of rod structure of a solution is discussed in Section 3; the main section is 4 where a catalogue of examples of $P$ is presented, some original and some collected from earlier work; finally there is a discussion of the inverse problem, the problem of obtaining $P$ directly from the rod structure and asymptotic quantities, in Section 5. Some extra material, including a brief consideration of recovering the metric from the patching matrix, is gathered in two appendices.

{\bf{Acknowledgements}}

I gratefully acknowledge assistance from Maciej Dunajski, Lionel Mason and Nick Woodhouse in preparing this review and catalogue, and Lars Andersson for the invitation to speak about this construction at a meeting at the Mittag-Leffler Institute in July 2023, the preparation for which prompted this review.

\section{Stationary, axisymmetric vacua and their twistor theory}

\subsection{Stationary axisymmetric vacua}

We start by reviewing the Einstein vacuum equations for stationary, axisymmetric, Lorentzian solutions, following \cite{ES}, \cite{MW} and \cite{W}. Suppose the two Killing vectors are $K=\partial/\partial t$ for time-translation and $L=\partial/\partial \phi$ for axisymmetry and that these commute. Also that the isometry group is \emph{orthogonally transitive}, that is the 2-surfaces orthogonal to the suraces of transitivity are integrable. For vacuum solutions the two scalars $*(K\wedge L\wedge dK)$ and  $*(K\wedge L\wedge dL)$ are necessarily constant and they will vanish if either $K$ or $L$ has a fixed point (which $L$ usually does have); then their common vanishing is the condition for orthogonal transitivity, which in turn ensures that coordinates can be chosen so that the metric is block $2\times 2$.

Now we can write the metric in the form

\be\label{1}g=h_{ij}dx^idx^j+
\left(\begin{array}{cc}
dt&d\phi\\
\end{array}\right)J
\left(\begin{array}{rr}
dt\\
d\phi\\
\end{array}\right)   \ee
with $i,j=1,2$ and $J$, sometimes called the \emph{Gram} or \emph{Yang} matrix, is
 \[J=\left(\begin{array}{cc}
        g_{tt} & g_{t\phi} \\
        g_{t\phi} & g_{\phi\phi}\\
\end{array}\right)=
\left(\begin{array}{cc}
        U & V \\
        V & W\\
\end{array}\right)\mbox{  say.}
\]
With Lorentzian signature $(+---)$, $h_{ij}$ is negative definite  and $J$ has negative determinant. We can include Riemannian metrics in the form (\ref{1}) by taking $h_{ij}$ and $J$ to be positive definite. In this case usually both Killing vectors have closed trajectories and the metric is said to be \emph{toric} (because of its torus symmetry).

It follows from the vacuum equations that det$J$ is harmonic with respect to the metric $h_{ij}$. We shall restrict to solutions with det$J$ nonconstant -- which could fail for example for solutions with two translation Killing vectors -- and then introduce the coordinate $r$ by
\be\label{2}
r^2=-\mbox{det}J,\mbox{   (Lorentzian)}\ee
or in the Riemannian case by
\be\label{3}
r^2=\mbox{det}J,\mbox{   (Riemannian)}.\ee
In either case, we choose $z$ to be the harmonic conjugate of $r$ (with respect to the 2-metric $h_{ij}$) and then the metric $h_{ij}$ becomes
\be\label{4}
h_{ij}dx^idx^j=\pm\Omega^2(dr^2+dz^2),\ee
with plus for Riemannian and minus for Lorentzian.  Note that $r$ is uniquely specified up to sign and $z$ up to sign and an additive constant. The locus of vanishing $r$, which we consider in detail in Section 3, is clearly of crucial interest and can be thought of for now as the axis of symmetry. The coordinates $(r,z)$ thus defined are often referred to as Weyl-Papapetrou coordinates.

It's quite common in the literature to write out (\ref{1}) at greater length in the Lorentzian case, in what is also called the Weyl-Papapetrou form, as
\be\label{5}
g=-f^{-1}(e^{2k}(dr^2+dz^2)+r^2d\phi^2)+f(dt-\omega d\phi)^2,\ee
when $\Omega^2=f^{-1}e^{2k}$ and then \cite{FW}, \cite{MW}
\be\label {6}
J=\left(\begin{array}{cc}
        f & -\omega f \\
        -\omega f & \omega^2f-\frac{r^2}{f}\\
\end{array}\right),\;\;\mbox{ (Lorentzian)}\ee
and $\mbox{det}J=-r^2$ so that the $r,z$ in (\ref{5}) are the Weyl-Papapetrou coordinates. The vacuum field equations for the metric (\ref{5}) can be found (with minor changes of  convention) in \cite{ES} (which can be taken as the reference for this part) as 
\bea\label{7}
\Delta(\log f)&=&-\frac{f^2}{r^2}(\omega_r^2+\omega_z^2)\\\label{8}
\left(\frac{f^2\omega_r}{r}\right)_r+\left(\frac{f^2\omega_z}{r}\right)_z&=&0\\\label{9}
k_r&=&\frac{r}{4f^2}(f_r^2-f_z^2)-\frac{f^2}{4r}(\omega_r^2-\omega_z^2)\\\label{10}
k_z&=&\frac{r}{2f^2}f_rf_z-\frac{f^2}{2r}\omega_r\omega_z
\eea
where
\[\Delta F(r,z)=F_{rr}+\frac1rF_r+F_{zz},\]
which is the 3-dimensional flat-space Laplacian in cylindrical polar coordinates for an axisymmetric function $F$. The integrability condition for (\ref{9}, \ref{10}) is satisfied by virtue of (\ref{7}, \ref{8})

Remarkably, (\ref{7}) and (\ref{8}) can be written as the matrix equation
\be\label{11}
r^{-1}\partial_r(rJ^{-1}\partial_rJ)+\partial_z(J^{-1}\partial_zJ)=0,\ee
\cite{LW}, usually called Yang's equation, and then the equations (\ref{9}, \ref{10}) can be written as a single (complex) equation for $\Omega$:
\be\label{12}
2i\partial_w\log(r\Omega^2)=r\mbox{Tr}(-J^{-1}J_wJ^{-1}J_w),\ee
where $w=z+ir$ and Tr means trace. As is to be expected, the integrability condition for (\ref{12}) is satisfied by virtue of (\ref{11}). It's worth noting that there is the freedom to add a constant of integration to $\log(r\Omega^2)$, which translates as a free multiplicative constant on $\Omega$. This can be thought of as the origin of the constant $k$ in (\ref{ct1}) below.

We'll be concentrating on (\ref{11}), but $\Omega$ and therefore (\ref{12}) are necessary for a discussion of conical singularities on the axis, which we'll be paying less attention to.

The next step towards solving the field equations has traditionally been, since \cite{Er}, to introduce a potential for $\omega$: (\ref{8}) is
 the integrability condition for the existence of a function $\psi$ satisfying
\[\psi_r=\frac{f^2}{r}\omega_z,\;\;\psi_z=-\frac{f^2}{r}\omega_r,\]
and fixed only up to an additive constant. These equations for $\psi$ turn out to be equivalent to
\[d\psi=*(K\wedge dK),\]
with $K=\partial_t$ as before. Thus $\psi$ is often referred to as the \emph{twist potential}  for the Killing vector $K$, since $*(K\wedge dK)$ can be called the \emph{twist} of the Killing vector, and this is co-closed by virtue of the vacuum equations. Evidently one could define a different twist potential by starting with a different Killing vector, and we'll take this up later (see (\ref{gh5})).  However, starting with $K$ we may consider the matrix \cite{FW}, \cite{MW}, \cite{LW}
\be\label{13}
J'=f^{-1}\left(\begin{array}{cc}
        1 & -\psi \\
        -\psi & \psi^2+f^2\\
\end{array}\right),\;\;\mbox{(Lorentzian)}\ee
then $J'$ also satisfies (\ref{11}) but now, note, with unit determinant.
The replacement of $J$ by $J'$ can be regarded as a B\"acklund transformation of the Yang equation (\ref{11}) (see \cite{MW} or \cite{NM1}). Note that one may recover $J$ from a knowledge of $J'$, and that there will be a variety of different $J'$ for the same $J$, obtained from different choices of Killing vector. 

It's worth noting that there is a version of (\ref{12}) involving only $J'$, namely
\be\label{12a}4ik_w=Tr(-J'^{-1}J'_wJ'^{-1}J'_w),\ee
with $k$ as in (\ref{5}).

 Following \cite{FW}, we shall call the matrix $J'$ the \emph{Ernst potential} for $J$ determined by $K$ (the more familiar scalar Ernst potential is the function $f+i\psi$, see e.g. \cite{ES}).

\medskip

So far we have considered Lorentzian metrics. The simplest way to shift to Riemannian metrics is to replace $J$ in (\ref{6}) by\footnote{This is done in \cite{AG} and seems to be due originally to Mason.}
\be\label{14}
J=\left(\begin{array}{cc}
        f & -\omega f \\
        -\omega f & \omega^2f+\frac{r^2}{f}\\
\end{array}\right),\;\;\mbox{(Riemannian)}\ee
when det$J$ becomes $r^2$, and $r,z$ are still Weyl-Papapetrou coordinates in the Riemannian metric
\be\label{15}g=f^{-1}(e^{2k}(dr^2+dz^2)+r^2d\phi^2)+f(dt-\omega d\phi)^2.\ee
Remarkably, (\ref{11}) and (\ref{12}) are both unchanged, but, with $\psi$ defined as before, $J'$ becomes
\be\label{16}
J'=f^{-1}\left(\begin{array}{cc}
        1 & -\psi \\
        -\psi & \psi^2-f^2\\
\end{array}\right),\;\;\mbox{(Riemannian)}\ee
now with determinant -1.

\subsection{Twistor construction of stationary axisymmetric vacuum solutions}

It was Witten's observation \cite{LW} that Yang's equation (\ref{11}) also determines anti-self-dual, static\footnote{In the literature of relativity `stationary' means possessing a time-like Killing vector, while `static' means possessing a hypersurface-orthogonal time-like Killing vector.}, axisymmetric gauge fields on $\mathbb{CM}$, and that these are solved by the Ward construction \cite{W2}, so we need a quick sketch of that. An anti-self-dual Yang-Mills field on the complexification of (a region of) Minkowski space $\mathbb{CM}$ defines a connection $\mathcal{A}$ which is flat on $\alpha$-surfaces i.e. on self-dual 2-planes. It therefore defines a vector space of constant sections on each $\alpha$-plane, and hence a holomorphic vector bundle $\mathcal{E}$ on the corresponding region in the space of $\alpha$-planes, which is projective twistor space $\mathbb{PT}=\mathbb{CP}^3$. Conversely the bundle $\mathcal{E}$ determines $\mathcal{A}$: restrict $\mathcal{E}$ to a projective line $\mathbb{CP}^1$ in $\mathbb{PT}$; it's part of standard twistor theory that these lines represent points in $\mathbb{CM}$ so by its construction $\mathcal{E}$ is trivial on these lines, and so this restriction splits; again it's part of standard twistor theory that the splitting determines $\mathcal{A}$ at the corresponding point in $\mathbb{CM}$. Reality conditions on $\mathcal{A}$ can be straightforwardly imposed.

If $\mathcal{A}$ is constant along a Killing vector $K$ of $\mathbb{CM}$, then $\mathcal{E}$ is constant along a corresponding holomorphic vector field on $\mathbb{CP}^3$ and one can contemplate taking quotients. This can similarly be done with two commuting symmetries, but now as Mason and Wodehouse \cite{MW} made clear, the action is not free and proper, so that the quotient, though a one-dimensional complex manifold, is not Hausdorff. However one may boldly go on and consider vector bundles on this non-Hausdorff Riemann surface. The pay-off is that, associated with a rod, that is to say a segment of the symmetry axis $r=0$,  on which the Killing vector $L$ vanishes and the Killing vector $K$ does not, the non-trivial patching matrix $P(z)$ of the bundle is obtained from the Ernst potential $J'(r,z)$ corresponding to $K$ simply by setting $r=0$:
\be\label{p1} P(z)=J'(0,z).\ee
It follows that there is a good deal of freedom in $P$, for example from changing the choice of $K$ or by adding a constant to $\psi$, which corresponds to conjugation of $P$ with a constant matrix.

A good deal more could be said here, but the original publications are clear and readily available.

\section{Rod structures}
The reference for this section is \cite{H}. 

\medskip

\noindent The locus $r=0$ in the $(r,z)$-half-plane from the previous section can loosely be regarded as the axis of symmetry, and transformed into a straight line (by the Riemann Mapping Theorem if necessary). It denotes the set of points where the rank of $J$ is strictly less than two, and contains a certain number of \emph{nodes} (also called \emph{nuts} in the literature, which risks confusion with nuts in the classification of \cite{GH1}, or \emph{turning points} - see \cite{CT}) where the rank is actually zero, connected by \emph{rods} on which the rank is one. The location of these nodes (in the $z$-coordinate), together with a 2-component vector label on each rod indicating which Killing vector lies in the kernel of $J$ there, is what is called the \emph{rod structure of the solution} \cite{H}. Together with asymptotic quantities, usually the mass and angular momentum, in the ALF case the rod structure characterises a solution (proved for five-dimensional AF Lorentzian in \cite{HY}; the case of four-dimensional AF Lorentzian is `black hole uniqueness' and in this context we can refer to \cite{kl}). Call this set of objects the \emph{data}, then there's no known prescription for generating all and only those data which are compatible with a solution. We'll give some examples.
\subsection{$\mathbb{E}^4$}
By this we mean flat Riemannian (i.e. Euclidean) four-space with the symmetry group chosen as two commuting rotations, in orthogonal planes (other choices are of course possible). Write the metric in terms of two sets of plane polars:
\be\label{e1}g=dR_1^2+dR_2^2+\left(\begin{array}{cc}
d\phi_1&d\phi_2\\
\end{array}\right)
\left(\begin{array}{cc}
        R_1^2 &0 \\
        0 & R_2^2\\
\end{array}\right)
\left(\begin{array}{rr}
d\phi_1\\
d\phi_2\\
\end{array}\right) ,\ee
then $r^2=\mbox{det}J=R_1^2R_2^2$ and so w.l.o.g. $r=R_1R_2$, when w.l.o.g. $z=\frac12(R_1^2-R_2^2)$. The axis of symmetry is the line $r=0$ and it has a node at the origin. The upper rod is at $R_2=0$, where $\partial_{\phi_2}$ vanishes and $z=\frac12R_1^2$ is positive; the lower rod is at $R_1=0$ where $\partial_{\phi_1}$ vanishes and $z=-\frac12R_2^2$ is negative. It's rather clear that  there could be no nonsingular solution with this rod structure and positive mass, which serves to show that some data are incompatible with the existence of a solution.

This example also typifies more complicated but \emph{asymptotically Euclidean} (AE) or \emph{asymptotically locally Euclidean} (ALE) solutions. Note the outermost rods are defined by a pair of linearly independent Killing vectors, which therefore automatically provide a basis of the Killing vectors.
\subsection{The Schwarzschild solution}
For the Lorentzian Schwarzschild metric in the usual coordinates, $J$ is given by
\[J=\left(\begin{array}{cc}
        1-2m/R &0 \\
        0 & -R^2\sin^2\theta\\
\end{array}\right)\]
and we'll take $m>0$. One speedily obtains
\[r=\sqrt{R(R-2m)}\sin\theta,\;z=(R-m)\cos\theta,\]
from which the rod structure can be read off by setting $r=0$: there are 3 rods, respectively at (i) $\theta=0$, (ii)$ R=2m$ and (iii) $\theta=\pi$, with nodes at $z=\pm m$. The central rod is the horizon, but in Riemannian Schwarzschild is the prime example of a \emph{bolt} in the classification of \cite{GH1}, being  a 2-sphere fixed by the isometry group.

For the top rod we have $\theta=0$ and so $z=R-m$, and since the Killing vector $\partial/\partial t$ is hypersurface-orthogonal we may take $\psi=0$. Thus the patching matrix $P_+$ from (\ref{13}) (and given already in  \cite{F}, \cite{FW} and \cite{W}) is

\[P_+=\left(\begin{array}{cc}
        \frac{z+m}{z-m} &0 \\
       0 & \frac{z-m}{z+m}\\
\end{array}\right).\]
We can use the algorithm from the next section for going past a node, to obtain the patching matrices for all rods as

\[P_0=\left(\begin{array}{cc}
        4(z^2-m^2) &0 \\
       0 & (4(z^2-m^2))^{-1}\\
\end{array}\right)\]
\[P_-=\left(\begin{array}{cc}
        \frac{z-m}{z+m} &0 \\
       0 & \frac{z+m}{z-m}\\
\end{array}\right).\]
For Riemannian Schwarzschild flip the sign on the lower right term in each $P_i$.

\subsection{Information from asymptotics}
The Schwarzschild solution, Lorentzian or Riemannian, is asymptotically flat and so has the asymptotics identified in \cite{H} for stationary, axisymmetric Ricci-flat metrics (though in the Lorentzian case these must have been written down earlier). In terms of a spherical radial coordinate $R$ these have the expansion (see also \cite{bs} or \cite{kl}) 
\[U=1-\frac{2m}{R}+O(R^{-2}),\;\;V=\frac{2L}{R}\sin^2\theta(1+O(R^{-1})),\;\;W=R^2\sin^2\theta(1+O(R^{-1})),\]
with $m,L$ the total mass and angular momentum (and this is now the Riemannian case). The Weyl-Papapetrou coordinates are
\[r=(R(R-2m))^{1/2}\sin\theta(1+O(R^{-1})),\;\;z=(R-m)\cos\theta(1+O(R^{-1})),\]
and the twist potential for the Killing vector $K=\partial_t$ is
\[\psi=\frac{2L}{R^2}\cos\theta(1+O(R^{-1})).\]
The outer rods are immediately identifiable as $\theta=0,\pi$, though we can have no information about nodes. For the top rod at $\theta=0$ we obtain
\be\label{h1}P_+=\left(\begin{array}{cc}
        1+\frac{2m}{z}+O(z^{-2}) &-\frac{2L}{z^2}+O(z^{-3}) \\
      - \frac{2L}{z^2}+O(z^{-3}) & -1+\frac{2m}{z}+O(z^{-2})\\
\end{array}\right).\ee
This expansion will be useful below. In particular, it shows that for Riemannian ALF we may suppose wlog that $P_+\rightarrow \left(\begin{array}{cc}
        1 & 0\\
        0 &-1\\
\end{array}\right)$ as $z\rightarrow\infty$. This was shown in \cite{F} for the Lorentzian case when instead $P_+\rightarrow I$.

\medskip

Having obtained the asymptotics for AF/ALF, it's worth taking a moment to do the corresponding thing for AE/ALE. Following a suggestion of James Lucietti, we can obtain the asymptotics of 4-dimensional AE/ALE from those implicit in \cite{H} for those particular 5-dimensional AF/ALF vacua which are just 4-dimensional vacua times a line. We consider the perturbation of (\ref{e1}), changing the coordinates to eliminate subscripts:
\be\label{e2}g=\Theta^2(du^2+dv^2)+\left(\begin{array}{cc}
d\phi&d\chi\\
\end{array}\right)
\left(\begin{array}{cc}
        u^2(1+\Phi) & \Xi\\
        \Xi & v^2(1-\Phi)\\
\end{array}\right)
\left(\begin{array}{rr}
d\phi\\
d\chi\\
\end{array}\right) ,\ee
where $\Phi,\Xi$ are to be regarded as small, and $\Theta$ is close to one (and we've used a gauge transformation to set the trace of the metric perturbation to zero). For the Weyl-Papapetrou coordinates we obtain
\[r=uv(1+h.o.),\;\;z=\frac12(u^2-v^2+h.o.).\]

Linearising the Yang equation (\ref{11}) with this $J$  leads to the Laplace equation for $\Phi$
\be\label{e3}
(uv\Phi_u)_u+(uv\Phi_v)_v=0,
\ee
and an equation for $\Xi$:
\be\label{e4}
(\frac{1}{uv}\Xi_u)_u+(\frac{1}{uv}\Xi_v)_v=0.\ee
This last is simplified by the substitution
\[S_v=\frac{1}{uv}\Xi_u,\;\;S_u=-\frac{1}{uv}\Xi_v,\]
when it becomes the Laplace equation for $S$:
\[(uvS_u)_u+(uvS_v)_v=0.\]
We choose a monopole solution for $\Phi$ and a symmetric dipole-like solution for $S$:
\be\label{e5}\Phi = \frac{2M}{(u^2+v^2)},\;\;S=\frac{L(u^2-v^2)}{(u^2+v^2)^3},\ee
when
\be\label{e7}\Xi=\frac{2Lu^2v^2}{(u^2+v^2)^3}.\ee
With these substituted into (\ref{e2}) we obtain the asymptotics given in \cite{H}, with $M,L$ here related to $\eta,\zeta$ there\footnote{It may be misleading to call these quantities $M$ and $L$ -- a more neutral usage would be $\eta,\zeta$ but, up to constant factors, they do occupy the places filled by $M,L$ in $P_+$ in the AF/ALF cases.}

For the corresponding $P$, we take the top rod to be $v=0$ so that the Killing vector $\partial_\chi$ vanishes there, and compute the twist potential of the Killing vector $\partial_\phi$ to be
\[\psi=-\frac{2Lu^4}{(u^2+v^2)^3}.   \]
On the axis where $v=0$, so $z=u^2/2$, we obtain
\be\label{e8}
P_+=\left(\begin{array}{cc}
        \frac{1}{2z}\left(1-\frac{M}{z}+h.o.\right) & \frac{L}{2z^2}+h.o.\\
         \frac{L}{2z^2}+h.o. & -2z(1+\frac{M}{z}+h.o.)\\
\end{array}\right),
\ee
which we'll use below.

\subsection{Patching matrices from rod structure and $J$}\label{sw3.3}
This is the starting point for the inverse problem. I'll quote a result from \cite{F} for an operation conveniently called `passing by a node'. Suppose then we have a rod structure with $k$ nodes located at $z=a_i,i=1,\ldots,k$ and a basis $(K,L)$ of Killing vectors. We know $J$ and on a rod on which $K$ does not vanish we can calculate $J'$. This will be singular at the next node, but it will at most have simple poles there (otherwise there would be a degenerate horizon in the Lorentzian case or a singularity in the Riemannian case) and the procedure for passing this node is as follows: suppose the node is at $z=0$ then there is freedom $P\rightarrow CPC^T$ with constant $C$ of unit determinant to set
\be\label{18}P=P^+= \left(\begin{array}{cc}
        (2z)^{-1}a(z) & b(z) \\
       b(z) & (2z)c(z)\\
\end{array}\right)\mbox{  with holomorphic }a,b,c\ee
above the node, when below the node we set
\[P=P^-= \left(\begin{array}{cc}
     (2z)a(z) & b(z) \\
 b(z) & (2z)^{-1}c(z)\\
\end{array}\right),\]
in words, the pole becomes a zero and the zero a pole.

From this prescription for passing by a node, we see that poles can only be introduced or removed at a node, so in general the entries in $P(z)$ can only be rational functions with simple poles at the nodes (poles off the axis give rise to singularities in the space-time metric -- see Appendix 2). If the nodes are located at $z=a_i$ for $i=1,\ldots,k$ then {\it{a priori}} the denominators of all entries in $P$ are $D(z)=\Pi_i^k(z-a_i)$ and numerators are polynomials, furthermore with their degrees fixed by the asymptotics. Also in AF/ALF or AE/ALE cases some coefficients in a series expansion are fixed by the asymptotic quantities. The determinant condition on $P(z)$ then fixes more, and in simple cases all, of the remaining coefficients. We'll return to this question in Section 5.

\section{Examples}
Now we present a menagerie of examples, paying particular attention to the Hermitian ones because of \cite{BG}.
\subsection{$\mathbb{E}^4$}
As we saw above, on the upper rod $z=\frac12R_1^2$ and the metric component $g_{\phi_1\phi_1}$ is $R_1^2=2z$ so from (\ref{16}) the patching matrix for the upper rod is
\be\label{ee1}P_+=\left(\begin{array}{cc}
        (2z)^{-1} & 0 \\
       0 & -(2z)\\
\end{array}\right).\ee
Following the same algorithm for $P_-$ gives
\[P_-=\left(\begin{array}{cc}
        -(2z)^{-1} & 0 \\
       0 & (2z)\\
\end{array}\right),\]
while taking $P_+$ past the node at $z=0$ would give
\[\widetilde{P}_-=\left(\begin{array}{cc}
        2z & 0 \\
       0 & -(2z)^{-1}\\
\end{array}\right),\]
which we regard as equivalent, since
\[P_-= \left(\begin{array}{cc}
        0 & 1 \\
       1 & 0\\
\end{array}\right)\widetilde{P}_-\left(\begin{array}{cc}
        0 & 1 \\
       1 & 0\\
\end{array}\right).\]
There are obviously no conicality issues in this case, and the explicit $P_+$ raises the expectation that generally for AE/ALE, $P_+$ does not have a limit as $z\rightarrow\infty$ but rather has terms with a simple zero and terms with a simple pole. We'll see an example like this with the C-metric in Section \ref{4.7}.

\subsection{Kerr}
The Lorentzian Kerr metric taken from \cite{HE} but with signature switched to $(+---)$  is
\be\label{k1}g_L=\left(1-\frac{2mR}{\Sigma^2}\right)dt^2+\frac{4maR\sin^2\theta}{\Sigma^2}d\phi dt-\left((R^2+a^2)\sin^2\theta+\frac{2mRa^2}{\Sigma^2}\sin^4\theta\right)d\phi^2-\Sigma^2\left(\frac{dR^2}{\Delta}+d\theta^2\right),\ee
in terms of two constants $m,a$ and with $\Sigma^2=R^2+a^2\cos^2\theta$ and $\Delta=R^2-2mR+a^2$, and is asymptotically-flat (AF). In the Lorentzian case, we'll require $m^2>a^2$ to avoid having a degenerate horizon, but we shall mainly be concerned with the Riemannian form of the metric, obtained by the substitution $(t,a)\rightarrow(it,ia)$ and a switch of overall sign (the Lorentzian case, up to the derivation of the patching matrices $P_i$, is given in \cite{F}, \cite{FW}). This is
\be\label{k2}
g_R=\left(1-\frac{2mR}{\Sigma^2}\right)dt^2+\frac{4maR\sin^2\theta}{\Sigma^2}d\phi dt+\left((R^2-a^2)\sin^2\theta-\frac{2mRa^2}{\Sigma^2}\sin^4\theta\right)d\phi^2+\Sigma^2\left(\frac{dR^2}{\Delta}+d\theta^2\right),\ee
where now $\Sigma^2=R^2-a^2\cos^2\theta$ and $\Delta=R^2-2mR-a^2$, and the metric is AF.

In the terminology of Section 2, we have
\be\label{k3} J= \left(\begin{array}{cc}
     1-\frac{2mR}{\Sigma^2} & \frac{2maR\sin^2\theta}{\Sigma^2} \\
 \frac{2maR\sin^2\theta}{\Sigma^2} &(R^2-a^2)\sin^2\theta-\frac{2ma^2R}{\Sigma^2}\sin^4\theta \\
\end{array}\right),\ee
and then
\be\label{k4}
r=(\mbox{det}J)^{1/2}=(R^2-2mR-a^2)^{1/2}\sin\theta,\ee
from which one rapidly finds
\be\label{k5}
z=(R-m)\cos\theta,\ee
having made a convenient choice of additive constant in $z$.

From (\ref{k4}) one rapidly identifies the rods: there are three, respectively (i) at $\theta=0$, (ii) at $R=R_+$ and (iii) at $\theta=\pi$, defining two nodes at $R=R_+$ and $\theta= 0$ or $\pi$, so also at $z=\pm(m^2+a^2)^{1/2}$. On the top and bottom rods, the Killing vector $L=\partial_\phi$ vanishes; on the middle rod the Killing vector $\partial_t-\frac{a}{2mR_+}\partial_\phi$ vanishes (this is in the Riemannian case when the middle rod is a bolt in the terminology of \cite{GH1}; in the Lorentzian case this is the location of the outermost Killing horizon where the Killing vector     $\partial_t+\frac{a}{2mR_+}\partial_\phi$ is null).

Adapting to the top rod, we'll calculate the twist potential for the Killing vector $K=\partial_t$. From the definition
\[d\psi=*(K\wedge dK),\]
we find
\be\label{k6}\psi=\frac{2ma}{R^2-a^2\cos^2\theta}.\ee
For the top rod, by evaluating on the axis, we find $R=z+m$ and
\be\label{k7}
P_+=J'(0,z)=f^{-1}\left(\begin{array}{cc}
        1 & -\psi \\
        -\psi & \psi^2-f^2\\
\end{array}\right)=
\frac{1}{(z^2-\sigma^2)}\left(\begin{array}{cc}
        (z+m)^2-a^2 & -2ma \\
        -2ma & -(z-m)^2+a^2\\
\end{array}\right)
\ee
where $\sigma=+(m^2+a^2)^{1/2}$ and the top node is at $z=\sigma$. 

In the corresponding way, we find on the lower rod
\be\label{k8}
P_-=J'(0,z)=
\frac{1}{(z^2-\sigma^2)}\left(\begin{array}{cc}
        (z-m)^2-a^2 & -2ma \\
        -2ma & -(z+m)^2+a^2\\
\end{array}\right).
\ee
Note that we could have obtained $P_-$ from $P_+$ by remarking that the Kerr metric has a discrete isometry $\theta\rightarrow -\theta$ or $z\rightarrow -z$.

For the middle rod, the matrix $P_+$ in (\ref{k7}) is not in the form of (\ref{18}) for passing by the top node but if we take
\[C=\left(\begin{array}{cc}
        1 & 0 \\
       \gamma & 1\\
\end{array}\right)\]
with $\gamma=a/(m+\sigma)$ then
\[P_+\rightarrow CP_+C^T=\left(\begin{array}{cc}
        \frac{(z+m)^2-a^2}{(z+\sigma)(z-\sigma)} & \frac{a(z+2m+\sigma)}{(m+\sigma)(z+\sigma)}\\
        \frac{a(z+2m+\sigma)}{(m+\sigma)(z+\sigma)} & \frac{-2m(z-\sigma)}{(m+\sigma)(z+\sigma)}\\
\end{array}\right),\]
which now is in the right form to apply (\ref{18}) and obtain $P_0$. 

The expression (\ref{k7}) for $P_+$, which may be compared with the expansion in (\ref{h1}), differs from the expression in \cite{F,FW,MW} as we are considering the Riemannian i.e. Wick-rotated, Kerr. For completeness we'll quote the $P_+$ matrix for {\bf{Lorentzian Kerr}} from these references:
\be\label{k10}
P_+=
\frac{1}{(z^2-\sigma^2)}\left(\begin{array}{cc}
        (z+m)^2+a^2 & 2ma \\
        2ma & (z-m)^2+a^2\\
\end{array}\right),
\ee
with $\sigma^2=m^2-a^2$, and recall this $a$ is $\pm i$ times the $a$ in (\ref{k7}).

\subsubsection{Hermiticity of Kerr}
Since both the self-dual and the anti-self-dual Weyl tensors are type D, the Riemannian Kerr metric has two distinct Hermitian structures, one associated with each orientation. The details can be found in \cite{PT}.

\subsubsection{Conicality issues}
There are known to be no conicality problems in the Lorentzian case: the central rod is the black-hole horizon and extension through this is well-known. For the Riemannian case see the discussion in \cite{GH1}.

\subsection{Taub-NUT}
This metric is discussed at some length in \cite{HE} and the elucidation of its structure was influential in the development of mathematical general relativity. It was first found as what would now be called a vacuum spatially homogeneous cosmology of Bianchi-type IX, \cite{T}. The metric form becomes singular where the surfaces of homogeneity become null but there exist extensions across these surfaces, where the character of the solution changes and it was seen to be identical with a solution independently found in \cite{NUT}. The standard name for the solution is obtained from these two references (\cite{NUT} gives NUT=Newman-Unti-Tamburino), but the term \emph{nut} in the classification of \cite{GH1} was motivated by this name. We shall take the Lorentzian form of the metric from \cite{HE}, with the signature switched to $(+---)$ and a slight notational change, as
\be\label{tn1}
g_L=\frac{dt^2}{A(t)}-A(t)(d\chi+2\ell\cos\theta d\phi)^2 -(t^2+\ell^2)(d\theta^2+\sin^2\theta d\phi^2),\ee
in terms of two constants $m,\ell$ and where $A(t)=-(t^2-2mt-\ell^2)(t^2+\ell^2)^{-1}$. Factorise the numerator of $A(t)$ as $(t_+-t)(t-t_-)$ with $t_\pm=m\pm(m^2+\ell^2)^{1/2}$ then for $t_-<t<t_+$ this is a Lorentzian Bianchi IX metric (the Taub solution) but the surfaces of homogeneity become null at $t=t_\pm$. The extensions through these surfaces, to the NUT part of the solution where $t$ becomes a radial coordinate, are discussed in  \cite{HE}, but we will be concerned with the Riemannian form of the metric which is obtained by the replacements $(t,\chi,\ell)\rightarrow (R,i\chi, iN)$ and a switch of overall sign (I'm using $N$ rather than $n$ to avoid confusion with a later reference to integer $n$). This is
\be\label{tn2}
g_R=\frac{dR^2}{U(R)}+U(R)(d\chi+2N\cos\theta d\phi)^2 +(R^2-N^2)(d\theta^2+\sin^2\theta d\phi^2),\ee
now with $U=(R^2-2mR+N^2)(R^2-N^2)^{-1}=-A$.
In the terminology of Section 2, we have
\be\label{tn3} J= \left(\begin{array}{cc}
     U & 2NU\cos\theta \\
 2NU\cos\theta & 4N^2U\cos^2\theta+(R^2-N^2)\sin^2\theta\\
\end{array}\right),\ee
and then
\be\label{tn4}
r=(\mbox{det}J)^{1/2}=(R^2-2mR+N^2)^{1/2}\sin\theta,\ee
from which one rapidly finds
\be\label{tn5}
z=(R-m)\cos\theta,\ee
with a convenient choice of additive constant in $z$.
 Now there are evidently three rods: call them respectively (i) $\theta=0$, (ii) $R=R_+=m+(m^2-N^2)^{1/2}$ and (iii) $\theta=\pi,$ and two nodes at respectively $(R,\theta)=(R_+,0)$ so $z=(m^2-N^2)^{1/2}$ and  $(R,\theta)=(R_+,\pi)$ so $z=-(m^2-N^2)^{1/2}$. We'll assume $m^2>N^2>0$ to keep the nodes real and separate, and it's convenient to write $\sigma=(m^2-N^2)^{1/2}$, so the nodes are at $z=\pm\sigma$. One can think of the parameter $m$ as a mass and then $N$ is usually called \emph{the nut parameter} (and is sometimes called the \emph{nuttiness}).

By inspection of (\ref{tn3}) we can see that the kernel of $J$ is $\partial_\phi-2N\partial_\chi$ on (i), $\partial_\chi$ on (ii) and $\partial_\phi+2N\partial_\chi$ on (iii). This should be compared with Kerr in the previous section: both have three rods and two nodes but the Killing vector labelling is quite different, so the rod structures are different, and in particular this metric won't be AF. In fact it is the paradigm ALF metric: consider surfaces of large constant $R$ in the metric (\ref{tn2}) then the first term vanishes, the last two terms give a round metric 2-sphere, and the middle term, with $U\rightarrow 1$, is the square of the connection one-form on a line bundle over the 2-sphere with curvature measured by $N$ -- this is ALF.

We'll calculate the twist potential for the Killing vector $K=\partial_\chi$. We find
\[d\psi=*(K\wedge dK)=2n\frac{(R^2-2mR+N^2)}{(R^2-N^2)^2}dR,\]
whence
\be\label{tn6}\psi=\frac{n(R-R_+)^2}{R_+(R^2-N^2)},\ee
where we've chosen the additive constant in $\psi$ so that $\psi(R_+)=0$. 
 For the top rod we set $\theta=0$ so $z=R-m$ and substituting into (\ref{16}) we obtain
\be\label{tn7}  P_+=J'(0,z)= \left(\begin{array}{cc}
     \frac{(z+m)^2-N^2}{(z-\sigma)(z+\sigma)} & -\frac{N(z-\sigma)}{R_+(z+\sigma)} \\
  -\frac{N(z-\sigma)}{R_+(z+\sigma)}  & -\frac{2\sigma(z-\sigma)}{R_+(z+\sigma)}\\
\end{array}\right).\ee
The patching matrices for the other rods are straightforward to obtain by the method of subsection \ref{sw3.3}.

Note that $P_+$ in (\ref{tn7}) tends to a constant as $z\rightarrow \infty$ but doesn't tend to the simple form $ \left(\begin{array}{cc}
     1 &0 \\
  0  & -1\\
\end{array}\right)$ that one expects for AF or ALF. We can impose this limit by replacing $P_+$ by $\tilde{P}_+=C^TP_+C$ with
\[C= \left(\begin{array}{cc}
     1 &\frac{N}{R_+} \\
  0  & 1\\
\end{array}\right).\]
This transformation is equivalent to adding a constant to $\psi$. It evidently preserves the symmetry and the determinant of $P_+$ and gives an equivalent twistor construction. We find
\be\label{tn77} \tilde{P}_+= \frac{1}{(z^2-\sigma^2)}\left(\begin{array}{cc}
      (z+m)^2-N^2 & 2Nz \\
  2Nz   & -(z-m)^2+N^2)\\
\end{array}\right).
\ee
This should be compared with (\ref{k8}): the nut-parameter term, in the off-diagonal slot,  here appears with a higher power of $z$ than the angular momentum term in (\ref{k8}), though it still goes to zero at large $z$ -- the asymptotic form of the patching matrix is as it was for Kerr.

\medskip

Since Taub-NUT also has a Lorentzian form, for completeness we'll give the $P_+$ for this: the $P_+$ matrix for {\bf{Lorentzian Taub-NUT}} is
\be\label{tn11} P_+= \frac{1}{(z^2-\sigma^2)}\left(\begin{array}{cc}
      (z+m)^2+\ell^2 & 2\ell z \\
  2\ell z   & (z-m)^2+\ell^2)\\
\end{array}\right),
\ee
where now $\sigma^2=m^2+\ell^2$ and $\ell^2=-N^2$ (as given in \cite{MW}).

\subsubsection{Self-dual Taub-NUT and Taub-bolt}

Self-dual Taub-NUT is the particular case $m=N$, so $R_+=m$ and $\sigma=0$: the Weyl tensor becomes self-dual and the two nodes apparently merge, but we don't get a second-order pole in $P_+$, rather from (\ref{tn7}) we find
\be\label{tn8}  P_+= \left(\begin{array}{cc}
     1+\frac{2m}{z} & -1 \\
  -1  & 0\\
\end{array}\right),\ee
so that there is a single node, and it is located at the origin (apparently resembling $\mathbb{E}^4$, but the asymptotics are quite different, as is $P$).
We'll see below that this form of $P_+$  is characteristic of the Gibbons-Hawking metrics, which include self-dual Taub-NUT. There are now only two rods and on the lower one we can directly calculate the patching matrix as
\be\label{tn9}  P_-= \left(\begin{array}{cc}
     1-\frac{2m}{z} & -1 \\
  -1  & 0 \\
\end{array}\right).\ee
Clearly we need a different rule from the one in section \ref{sw3.3} for going past a node with the Gibbons-Hawking metrics. We'll return to this in section \ref{4.5}.

\medskip

The so-called Taub-bolt metric was spotted as another particular case of the general Taub-NUT in \cite{Pa}. It is neither self-dual nor anti-self-dual, but is Hermitian as we'll see below. To obtain it, one sets $m=5N/4$, so $R_+=2N$ and it's convenient to work with $\sigma$ which is now $3N/4$. We find

\be\label{tn10} P_+= \left(\begin{array}{cc}
     \frac{(z+3\sigma)(z+\frac13\sigma)}{(z-\sigma)(z+\sigma)} & -\frac{(z-\sigma)}{4(z+\sigma)} \\
  -\frac{(z-\sigma)}{4(z+\sigma)}  & -\frac{3(z-\sigma)}{4(z+\sigma)}\\
\end{array}\right).\ee
In this example the central rod is a {\it bolt} in the terminology of \cite{GH1}, that is to say it corresponds to an $S^2$ in the manifold on which the Killing vector $K=\partial_\chi$ vanishes. The central rod of the Riemannian Kerr metric is also a bolt but the central node of the self-dual Taub-NUT solution is a {\it nut}, as it is an isolated zero of a Killing vector (in fact of $\partial_\chi$).

\subsubsection{Hermiticity of Taub-NUT}
Both SD and anti-SD Weyl tensors are type D for Taub-NUT, so it has two Hermitian structures, with opposite orientations. With the metric as in (\ref{tn2}) consider the basis of 1-forms
\[\Theta^0=\frac{dR}{\sqrt{U}},\;\Theta^1=(R^2-N^2)^{1/2}d\theta,\;\Theta^2=(R^2-N^2)^{1/2}\sin\theta d\phi,\;\Theta^3=\sqrt{U}(d\chi+2N\cos\theta d\phi),\]
then an integrable complex structure is defined by taking the holomorphic 1-forms to be
\[\Theta^0+i\Theta^3,\;\;\Theta^1+i\Theta^2.\]
The candidate K\"ahler form is
\[\Omega=\Theta^0\wedge\Theta^3+\Theta^1\wedge\Theta^2\]
and is not closed but
\[\Omega_K=(R-N)^{-2}\Omega\]
is, so that $(R-N)^{-2}g_R$ with $g_R$ as in (\ref{tn2}) is K\"ahler.

For the other orientation, switch the sign on $\chi$ and $N$ then the same choices give the other complex structure and now $(R+N)^{-2}g_R$ is K\"ahler.

\subsubsection{Conicality}
See \cite{GH1} for a discussion of conicality in self-dual Taub-NUT and the Eguchi-Hanson metric (coming up in Section \ref{4.5}); the absence of conical singularities in Taub-bolt is shown in \cite{Pa}.
\subsection{Kerr-Taub-bolt}
The Kerr-Taub-bolt metric of \cite{GP} is in this family, in the sense of having two nodes. Its Weyl tensor is type D, so it's Hermitian, and it has a certain amount of fall-off but it isn't in the list of regular AF Hermitian metrics in \cite{BG} so it must not be regular -- as also pointed out in \cite{CT} it has conical singularities.

Here it is anyway
\[g=X\left(\frac{dR^2}{\Delta}+d\theta^2\right)+\frac{\sin^2\theta}{X}(adt+P_Rd\phi)^2+\frac{\Delta}{X}(dt+P_\theta d\phi)^2,\]
with
\[X=R^2-(a\cos\theta+N)^2\]
\[P_R=R^2-a^2-\frac{N^2}{N^2-a^2}\]
\[P_\theta=-a\sin^2\theta+2n\cos\theta-\frac{an^2}{n^2-a^2}\]
\[\Delta=R^2-2mR+N^2-a^2.\]
Following the algorithm, calculate
\[r=\Delta^{1/2}\sin\theta,\;z=(R-m)\cos\theta\]
and then for the Killing vector $K=\partial_t$ the twist potential is
\[\psi=\frac{2}{X}(NR-m(a\cos\theta +N)).\]
Restrict to the upper axis $\theta=0$ and compute from (\ref{16})\footnote{Previously given in \cite{AG}.}:

\be\label{kn1}  P_+= \left(\begin{array}{cc}
     \frac{z^2+2mz+m^2-(a+N)^2}{z^2-\sigma^2} & \frac{2(Nz-am)}{z^2-\sigma^2} \\
   \frac{2(Nz-am)}{z^2-\sigma^2}  & \frac{-z^2+2mz-m^2+(a-N)^2}{z^2-\sigma^2}\\
\end{array}\right),\ee
where $\sigma^2=m^2+a^2-N^2$ and the nodes are at $\pm\sigma$. 

\medskip

The interest of this example is that it can also be obtained as follows: make a general two-node ansatz with asymptotic quantities $(m,a,N)$ for mass, angular momentum and nut parameter, and free parameters $(A,B,\sigma)$, as follows
\be\label{kn2}  P_+= \left(\begin{array}{cc}
     \frac{z^2+2mz+A}{z^2-\sigma^2} & \frac{2(Nz-am)}{z^2-\sigma^2} \\
   \frac{2(Nz-am)}{z^2-\sigma^2}  & \frac{-z^2+2mz-m^2+B}{z^2-\sigma^2}\\
\end{array}\right),\ee
now impose the determinant condition on $P_+$, to find that the remaining constants, including the location of the nodes, are fixed by the asymptotic quantities.

\subsection{The Gibbons-Hawking metrics}\label{4.5}
These are characterised as being hyper-K\"ahler with a triholomorphic Killing vector, or alternatively as Ricci-flat with a self-dual Weyl tensor and a Killing vector whose derivative is self-dual \cite{GH2, WT}. The metric can be written
\be\label{gh1} g=V(dx^2+dy^2+dz^2)+V^{-1}(dt+{\bf{\omega}}\cdot{\bf{dx}})^2,\ee
in terms of a function $V$ and a 3-vector ${\bf{\omega}}$, and the triholomorphic Killing vector is $K=\partial_t$. The three K\"ahler forms can be taken to be
\[\Omega_1=(dt+{\bf{\omega}}\cdot{\bf{dx}})\wedge dx-Vdy\wedge dz\]
and its two cyclic permutations, and then closed-ness of these implies the field equation
\[\partial_1V=\partial_2\omega_{3}-\partial_3\omega_{2}\mbox{ and cyclic permutations,}\]
or in a 3-vector notation
\[\nabla V=\nabla\times {\bf{\omega}}.\]
It follows that $V$ is harmonic w.r.t. the flat 3-metric.

It's worth noting that for some authors (e.g. \cite{CV}) a gravitational instanton is necessarily hyper-K\"ahler.

We are interested in toric symmetry so we need a second Killing vector which commutes with $K=\partial_t$. We don't a priori know what it has to be, so write it in its most general form as 
\[L=A\partial_x+B\partial_y+C\partial_z+D\partial_t,\]
and seek to restrict it. First off, commutation with $K$ at once forces $A,B,C,D$ to be independent of $t$, and $LV$ to be zero. Now $L$ must Lie-drag the tensor $V^{-1}g_{ab}-K_aK_b$, which is the flat 3-metric in (\ref{gh1});  $D$ must be a constant, so we can set it to zero by subtracting a multiple of $K$ from $L$, and then $L$ is just a symmetry of the flat 3-metric. It could be a translation, but then det$J$ would be a constant and there won't be a $J'$; it could be a screw-rotation, but that is incompatible with asymptotic flatness; so we'll take it to be a rotation, $L=\partial_\phi$. Adapting (\ref{gh2}) to this we see that ${\bf{\omega \cdot dx}}$ can be taken to be $Wd\phi$ when the metric is
\be\label{gh2}g=V(dr^2+r^2d\phi^2+dz^2)+V^{-1}(dt+Wd\phi)^2,\ee
with $V(r,z)$ and $W(r,z)$ related by
\[rV_z=W_r,\;\;rV_r=-W_z.\]
For $J$ we find
\[J=\left(\begin{array}{cc}
         V^{-1}& V^{-1}W\\
        V^{-1}W& V^{-1}W^2+r^2V\\
\end{array}\right),\]
so det$J=r^2$ and $r$ is already the Weyl-Papapetrou coordinate. There would seem to be only one rod, where $\partial_\phi$ vanishes, but this is an error: $V$ can have simple poles which are removable singularities of $J$ and are genuinely isolated zeroes of the Killing vector $\partial_t$, {\it nuts} in the terminology of \cite{GH1}. Now each of these nuts is a {\it node} in our earlier usage.

Following the method of Section 2, we can find the twist potential for $K$ as
\[\psi=\frac{1}{V}+\mbox{ constant}.\]
While it may sometimes be convenient to choose the constant to set $\psi=0$ at infinity, if we take the constant to be zero then from (\ref{16}) we obtain the simple expression
\be\label{gh3}
P=J'(0,z)=
\left(\begin{array}{cc}
         V(0,z)& -1\\
        -1& 0\\
\end{array}\right),\ee
and this should hold on all rods (compare (\ref{tn8})). Note that we only obtain this simple form because we are using the triholomorphic Killing vector.

\medskip

 The two classes of interesting (smooth and ALF or ALE) examples are the  multi-Taub-NUT and the multi-Eguchi-Hanson metrics\footnote{These two classes are the only regular GH solutions with AF/ALF/AE/ALE asymptotics: see \cite{Mi} and references therein.}. Both have a $V$ corresponding to a series of point masses all of mass $m$ at a sequence of nodes $z=a_i,i=1,\ldots,n$ along the $z$-axis, but with $V=1$ at large distances for multi-Taub-NUT and $V=0$ at large distances for multi-Eguchi-Hanson:
\bea V_{mTN}&=&1+\sum_{i=1}^n\frac{m}{(r^2+(z-a_i)^2)^{1/2}}\label{gh4}\\
V_{mEH}&=&\sum_{i=1}^n\frac{m}{(r^2+(z-a_i)^2)^{1/2}}.\label{gh5}
\eea
The top-left entry in $P$ for multi-Taub-NUT is now
\[V(0,z)=1+\sum_{i=1}^n\frac{m}{(z-a_i)},\]
and the same without the one for multi-Eguchi-Hanson, and as we proceed down the axis from large positive $z$, as we pass the node $a_i$ say we simply switch the term $m(z-a_i)^{-1}$ to $m(a_i-z)^{-1}$ in this sum.

{\bf Comment}: Note that, at a given integer $n$,{\em  mTN} and {\em mEH} have the same arrangement of rods, but different asymptotics and consequently different $P$.

\medskip

\subsubsection{Conicality}
All the angular singularities can be removed (see \cite{GH2}) so these metrics provide regular examples with as many rods as desired, and examples of both ALE (mEH) and ALF (mTN).

\subsubsection{Flat space another way}
 It's well-known that the Eguchi-Hanson metric in Gibbons-Hawking form but with $n=1$ is actually flat. In this case we can take
\be\label{gh7}
P_+=\left(\begin{array}{cc}
         1/z& -1\\
        -1& 0\\
\end{array}\right),\ee
which is different from (\ref{ee1}): the rod structure here, as in Section 3.1, has a single node, but here the Killing vectors vanishing on the upper and lower rods are the same, while in Section 3.1 they were different. Now consider the transformation
\be\label{gh6}P_+\rightarrow CP_+C^T\mbox{  with  } C=\left(\begin{array}{cc}
         \alpha^{-1}& 0\\
        0& \alpha\\
\end{array}\right)\left(\begin{array}{cc}
         1& 0\\
        \beta z/\alpha& 1\\
\end{array}\right)\left(\begin{array}{cc}
         1& 0\\
        \gamma/\alpha& 1\\
\end{array}\right)=\left(\begin{array}{cc}
         \alpha^{-1}& 0\\
        \beta z+\gamma& \alpha\\
\end{array}\right).\ee
Note all three factors of $C$ have determinant one, so det $P_+$ is unchanged. The first factor rescales the diagonal entries; the second effects a change in the choice of Killing vector on the upper rod\footnote{For the justification of this claim, see the discussion around eqn(9) in \cite{FW}.}; and the third corresponds to adding a constant to the twist potential $\psi$. The effect on $P_+$ is then
\[P_+\rightarrow 
\left(\begin{array}{cc}
         1/\alpha^2z& ((\alpha+\beta)z+\gamma)/(\alpha z)\\
        ((\alpha+\beta)z+\gamma)/(\alpha z)& ((2\alpha\beta+\beta^2)z^2+2\gamma(\alpha+\beta)z+\gamma^2)/z\\
\end{array}\right).\]
If we now take $\gamma=0,\alpha=-\beta=\sqrt{2}$ then $P_+$ becomes the $P_+$ of (\ref{ee1}). What is new here is the middle factor in (\ref{gh5}), and this is taken from \cite{FW}: the result there is that this is, quite generally, the way to modify $P$ under a change of choice of Killing vector. We'll make use of it in the discussion of Pleba\'nski-Demia\'nski metrics below.

One might prefer (\ref{gh7}) to (\ref{ee1}) as a patching matrix for flat space, and therefore as the limiting $P_+$ for an ALE metric, as it is bounded at large $z$, but it seems to make little difference in practice.

\subsubsection{Another particular Gibbons-Hawking metric}
This metric was discussed in \cite{AG} but originated in \cite{B}. The metric is
\be\label{ba1}
g=f(d\rho^2+d\theta^2)+f^{-1}\sin^2\theta(ad\chi+\sigma\cosh\rho d\phi)^2+f^{-1}\sinh^2\rho(bd\chi+\sigma\cos\theta d\phi)^2,\ee
with
\[f=b\cosh\rho -a\cos\theta,\]
and I've introduced a parameter $\sigma$ not present in \cite{B} (or perhaps set equal to one there, though it is useful to have it present for our purposes).
The fact that the metric is in this section implies that the metric is of Gibbons-Hawking type, which wasn't mentioned in \cite{B} and may not have been known. It was pointed out to me by Maciej Dunajski.

One calculates $J$ and then its determinant to find the Weyl-Papapetrou coordinates as
\[r= \sigma\sin\theta\sinh\rho,\;\;z=\sigma\cos\theta\cosh\rho,\]
and these are confocal-elliptic-hyperbolic coordinates: the curves of constant $\rho$ are ellipses, the curves of constant $\theta$ are hyperbolas, all with foci at $r=0,z=\pm \sigma$. Also
\[z+ir=\sigma\cosh(\rho+i\theta),\]
so that
\[dz^2+dr^2=\sigma^2(\cosh^2\rho-\cos^2\theta)(d\rho^2+d\theta^2).\]
Comparing (\ref{ba1}) with (\ref{gh1}) we conclude that
\[V=\frac{f}{\sigma^2(\cosh^2\rho-\cos^2\theta)}=\frac{(\alpha+\beta)\cosh\rho-(\beta-\alpha)\cos\theta}{\sigma^2(\cosh^2\rho-\cos^2\theta)}\]
where we've put $a=\beta-\alpha,b=\alpha+\beta$, so that
\be\label{ba2}V=\frac{\alpha}{\sigma(\cosh\rho-\cos\theta)}+\frac{\beta}{\sigma(\cosh\rho+\cos\theta)}.\ee

Should we need it, note that the other ingredient in (\ref{gh1}) is now
\[\omega=Ad\phi\mbox{   with   }A=\frac{\alpha(\cosh\rho\cos\theta-1)}{(\cosh\rho-\cos\theta)}+\frac{\beta(\cosh\rho\cos\theta+1)}{(\cosh\rho+\cos\theta)}.\]

Putting the expression for $V$  back in terms of $r,z$ we obtain
\[V=\frac{\alpha}{(r^2+(z-\sigma)^2)^{1/2}}+\frac{\beta}{(r^2+(z+\sigma)^2)^{1/2}}\]
which is the Newtonian potential for point masses $\alpha,\beta$ on the axis at $z=\sigma$ or $-\sigma$ respectively. This is like the Eguchi-Hanson metric with two nodes but with unequal masses. The metric will not be smooth, as already observed in \cite{B}, but the corresponding $P_+$ will be interesting below. From (\ref{gh3}) this is
\be\label{ba3}P_+=\left(\begin{array}{cc}
       \frac{bz-a\sigma}{z^2-\sigma^2}  & -1\\
        -1& 0\\
\end{array}\right) ,  \ee
which we return to in Section 5.

\subsection{The Weyl solutions and double Schwarzschild}
If we put $\omega=0$ in (\ref{5},{7},{13}) then also $\psi=0$ and these are the (Lorentzian) Weyl static solutions. Note that for these $f=e^{2U}$ with $U$ harmonic. As is well-known, the Schwarzschild solution lies in this family with $U$ being the potential for a notional massive rod (which is also a rod in the sense of rod structure) with mass $1/2$ per unit length. It is possible to superimpose any number of distinct but collinear such massive rods to produce a multi-Schwarzschild metric, though one cannot eliminate all the conical singularities, \cite{IK}.

In the simplest case, for the double Schwarzschild, suppose there are two massive rods lying along the $z$-axis, one of length $b_1$ with ends at $a_0$ and $a_1=a_0+b_1$ and the other of length $b_2$ with ends at $a_2$ and $a_3=a_2+b_2$, and we'll suppose the $a_i$ are in order so that $a_0<a_1<a_2<a_3$. In terms of rod structure, there are five rods and four nodes at $z=a_i$. Where the massive rods are, the Killing vector $\partial_t$ is null and on the other three rods $\partial_\phi$ vanishes. We'll label the rods from the bottom as $P_-,P_1,P_2,P_3,P_+$.

Following \cite{IK}, introduce
\[z_i=z-a_i,\;\rho_i=(r^2+z_i^2)^{1/2},\]
then consider
\be\label{w0}U=\frac12\log\left(\frac{\rho_0+\rho_1-b_1}{\rho_0+\rho_1+b_1}\right) +\frac12\log\left(\frac{\rho_2+\rho_3-b_2}{\rho_2+\rho_3+b_2}\right) .\ee
This is harmonic, and is the $U$ for double Schwarzschild from \cite{IK}. These authors also give an explicit expression for $k$, which we don't need.

At once
\[f=e^{2U}=\frac{(\rho_0+\rho_1-b_1)(\rho_2+\rho_3-b_2)}{(\rho_0+\rho_1+b_1)(\rho_2+\rho_3+b_2)},\]
and we can read off the patching matrices from (\ref{13}): on the axis with $z>a_3$ we have
\[\rho_0=z-a_0,\;\rho_1=z-a_1,\;\rho_2=z-a_2,\;\rho_3=z-a_3,\]
so that
\be\label{w1}f=\frac{(z-a_1)(z-a_3)}{(z-a_0)(z-a_2)},\ee
and
\be\label{w2}
P_+=\left(\begin{array}{cc}
         \frac{(z-a_0)(z-a_2)}{(z-a_1)(z-a_3)}& 0\\
        0& \frac{(z-a_1)(z-a_3)}{(z-a_0)(z-a_2)}\\
\end{array}\right).\ee
In the same way, on the middle rod where $a_1<z<a_2$ we arrive at
\[P_2=\left(\begin{array}{cc}
         \frac{(z-a_0)(z-a_3)}{(z-a_1)(z-a_2)}& 0\\
        0& \frac{(z-a_1)(z-a_2)}{(z-a_0)(z-a_3)}\\
\end{array}\right),\]
and on the bottom rod, where $z<a_0$
\[P_-=\left(\begin{array}{cc}
         \frac{(z-a_1)(z-a_3)}{(z-a_0)(z-a_2)}& 0\\
        0& \frac{(z-a_0)(z-a_2)}{(z-a_1)(z-a_3)}\\
\end{array}\right).\]
To obtain $P_1,P_3$ we follow the rule in section \ref{sw3.3} for going round a node. Taking $P_+$ round the node at $a_3$ we obtain
\[P_3=\left(\begin{array}{cc}
         \frac{4(z-a_0)(z-a_2)(z-a_3)}{(z-a_1)}& 0\\
        0& \frac{(z-a_1)}{4(z-a_0)(z-a_2)(z-a_3)}\\
\end{array}\right),\]
and taking $P_2$ round the node at $a_1$ we obtain
\[P_1=\left(\begin{array}{cc}
         \frac{4(z-a_0)(z-a_1)(z-a_3)}{(z-a_2)}& 0\\
        0& \frac{(z-a_2)}{4(z-a_0)(z-a_1)(z-a_3)}\\
\end{array}\right).\]
Having obtained the complete set of patching matrices, one can check that passing round the nodes at $a_0$ and $a_2$ is consistent.

\medskip

It's clear that one can obtain an arbitrary even number of nodes this way following the pattern of (\ref{w1},\ref{w2}) but, as we see next, there will always be conical singularities. This example is interesting as it is not Hermitian, equivalently neither of its Weyl spinors is type D.
\subsubsection{Conicality}
It's well-known that, while the horizons are regular (in the Lorentzian case), the conical singularities at the vanishing of $\partial_\phi$ cannot all  be removed. In the double Schwarzschild case, the periodicity of $\phi$ can be chosen to make the two outer rods nonsingular, when the central rod remains singular and is interpreted as a strut keeping the two black holes apart, or alternatively the periodicity of $\phi$ can be chosen to remove the conical singularity on the central rod when the two outer rods remain singular and are interpreted as ropes holding the two black holes apart. 

This discussion was for the Lorentzian metric, but it holds for Riemannian with the usual sign changes.

\subsection{C-metric}\label{4.7}
The C-metric (see e.g. \cite{ES}) was considered from the current point of view  in \cite{F} where it was seen to have three nodes. The metric is vacuum type D and so admits a valence-two Killing spinor: when made Riemannian, it will be Hermitian and conformal to K\"ahler in two ways. Following \cite{F} we can take the Lorentzian form of the metric to be
\[g_L=\frac{1}{A^2(x+y)^2}\left(F(y)dt^2-G(x)d\phi^2-\frac{dx^2}{G(x)}-\frac{dy^2}{F(y)}\right),\]
with $G(x)=1-x^2-\alpha x^3,F(y)=-1+y^2-\alpha y^3$ (so $F(y)=-G(-y)$), and real constants $\alpha,A$ where $m=A/(2\alpha)$ is the mass in the Lorentzian version. We'll assume that $m$ and $\alpha$ are both nonzero.

We obtain the Riemannian version by replacing $t$ by $it$ and switching the overall sign:
\be\label{c1}
g_R=\frac{1}{A^2(x+y)^2}\left(F(y)dt^2+G(x)d\phi^2+\frac{dx^2}{G(x)}+\frac{dy^2}{F(y)}\right).\ee
There are two distinct Hermitian structures defined by taking either $(idt+dy/F,id\phi+dx/G)$ or $(idt+dy/F,id\phi-dx/G)$ as holomorphic one-forms.
We note that the 2-forms
\[\omega_{\pm}=dt\wedge dy\pm d\phi\wedge dx\]
are K\"ahler w.r.t. one or other of the Hermitian structures for the rescaled metric $(x+y)^2g_R$ so that also the coordinates $x$ and $y$ are  Hamiltonians for the Killing vectors $\partial_\phi,\partial_t$ respectively.

The metric functions $F$ and $G$ both have 3 real roots if $\alpha^2<4/27$, which we'll assume. Label the roots of $G$ as $\beta_1<\beta_2<0<\beta$, then the roots of $F$ are $-\beta<0<-\beta_2<-\beta_1$ (also we can write $G(x)=-\alpha(x-\beta)(x-\beta_1)(x-\beta_2)$ or $F=-\alpha(y+\beta)(y+\beta_1)(y+\beta_2)$ whenever desired). For a Riemannian metric we need $F(y)$ and $G(x)$ non-negative and there are four regions in the $(x,y)$-plane that would satisfy this. However, for one of them the region doesn't meet the line $x+y=0$, which represents infinity (as is evident from (\ref{c1})), and in two more of them $x$ or $y$ can go to infinity (and the metric is singular), so we'll choose the fourth, which is the region $R$ defined by $\beta_2\leq x\leq \beta$ and $-\beta_2\leq y\leq-\beta_1$. Note that $F=0$ on two sides of this region and $G=0$ on the other two: the coordinate $r$ will vanish all the way round the rectangle. 

The Weyl-Papapetrou coordinates are \cite{F}
\[r=\frac{(F(y)G(x))^{1/2}}{A^2(x+y)^2},\;\;z=\frac{(\alpha xy^2 -\alpha x^2y -2xy-2)}{2A^2(x+y)^2}.\]
If we confine $(x,y)$ to the rectangle $\beta_2<x<\beta,-\beta_2<y<-\beta_1$, then the corner $(\beta_2,-\beta_2)$ is at infinity and the other three corners are nodes, so the four sides are the rods, and we see at once which Killing vector vanishes on each rod: $\partial_t$ when $F=0$ and $\partial_\phi$ when $G=0$.

Following \cite{F} we can calculate the $z$-coordinates of the three nodes to find
\[z(\beta,-\beta_2)=\mu:=\frac{\alpha\beta_1}{2A^2},\; z(\beta,-\beta_1)=\nu:=\frac{\alpha\beta_2}{2A^2},\; z(\beta_2,-\beta_1)=\xi:=\frac{\alpha\beta}{2A^2}  \]
with $\mu<\nu<0<\xi$. Since the Killing vectors are both twist-free, all the patching matrices will be diagonal (see (\ref{16})) and we just need to restrict $g_{tt}$ or $g_{\phi\phi}$ to the axis. On the top-most rod $z>\xi$ the Killing vector $\partial_\phi$ vanishes so we use $\partial_t$ to construct the patching matrix $P_+$ according to (\ref{16}), and we simply need
\[f=g_{tt}=\frac{F(y)}{A^2(\beta_2+y)^2}=\frac{-\alpha(y+\beta)(y+\beta_1)}{A^2(y+\beta_2)}\]
where we have substituted the $x$-value for this rod, namely $\beta_2$, and used the factorisation of $F(y)$. We want this in terms of $z$ where 
\[z=\frac{(\alpha \beta_2y^2 -\alpha \beta_2^2y -2\beta_2y-2)}{2A^2(\beta_2+y)^2}=\frac{(\alpha\beta_2^2y-2)}{2A^2\beta_2(y+\beta_2)}.\]
Eliminating $y$ from $f$ in favour of $z$ we find the simple expression
\be\label{c2}f=\frac{2(z-\mu)(z-\xi)}{(z-\nu)}.\ee
For the other three rods we employ the switching procedure of subsection \ref{sw3.3}. The function $f$ above is replaced in turn by
\[\frac{(z-\mu)}{2(z-\xi)(z-\nu)},\;\frac{2(z-\mu)(z-\nu)}{(z-\xi)},\;\mbox{  and  }\frac{(z-\nu)}{2(z-\mu)(z-\xi)},\]
as one passes the nodes at $\xi,\nu,\mu$ respectively. By the end, on the bottom rod, we've arrived at the inverse $P$ of the top rod. This can be compared with $\mathbb{E}^4$ in subsection 4.1 and contrasted with the Kerr case in subsection 4.2 and the double Schwarzschild in subsection 4.4 where the patching matrix $P_+$ for the top rod is the same as the patching matrix $P_-$ for the bottom rod. This is because the asymptotics are different: Kerr and double Schwarzschild are AF while this C-metric is AE. To see this for the Riemannian C-metric, we look near to the point $(x,y)=(\beta_2,-\beta_2)$: set
\[x=\beta_2+\frac{c}{R^2}\cos^2(\theta/2),\;\;y=-\beta_2+\frac{c}{R^2}\sin^2(\theta/2),\]
where $c$ is a constant to be chosen, and we're thinking of $R$ as large. Then
\[G(x)=k\frac{c}{R^2}\cos^2(\theta/2),\;F(y)=k\frac{c}{R^2}\sin^2(\theta/2)\]
where $k=G'(\beta_2)=F'(-\beta_2)$ which can be seen from the graph of $G$ to be positive, and we neglect higher order in $1/R$. Substitute for $F$ and $G$ into (\ref{c1}), along with $t=\alpha-\beta,\phi=\alpha+\beta$ to find
\[g_R=\frac{4k}{A^2c}\left(dR^2+\frac14R^2(d\theta^2+\sin^2\theta d\beta^2+(d\alpha+\cos\theta d\beta)^2))\right)\]
which is the standard flat metric on $\mathbb{R}^4$ if we set $c=4k/A^2$. 

Being static and axisymmetric, the C-metric is in Weyl's class and arises from a harmonic potential $U$ similar in form to that in (\ref{w0}) but with an odd number of summands. Following the pattern of the double-Schwarzschild, we may write down a $P$ with any odd number of nodes, for example for nodes at $z=a_1,\ldots,a_5$ take
\[f:=\frac{2(z-a_1)(z-a_3)(z-a_5)}{(z-a_2)(z-a_4)},\]
in place of (\ref{w1}) and (\ref{c2}). Again it will not be possible to remove all the conical singularities, but the metric will this time be ALE and still not Hermitian.

\medskip

There is a generalisation of the C-metric known as the spinning or twisting C-metric (see e.g. \cite{PD}). The metric is not diagonal and so neither is $P$ but this metric is a special case of the Pleba\'nski-Demia\'nski metrics which in turn are a special case of the Chen-Teo metric so we'll leave it until we look at those.
\subsubsection{Conicality}
The conical singularities on the axis cannot all be resolved (a familiar fact in the Lorentzian case where they are sometimes interpreted as struts holding two black holes apart). To see this, consider the rod at $x=\beta$ and calculate the periodicity in $\phi$ needed to make this regular: nearby $x=\beta-\rho^2$ so $G(x)=-G'(\beta)\rho^2+h.o.$ and the condition is
\[\Delta\phi=4\pi/(-G'(\beta)),\]
while for the rod at $x=\beta_2$ one needs
\[\Delta\phi=4\pi/(G'(\beta_2)),\]
and these cannot be the same:
\[G'(\beta)+G'(\beta_2)=-\alpha(\beta-\beta_2)^2\neq 0.\]
Similarly the rods at $y=-\beta_1,-\beta_2$ require
\[\Delta t =4\pi/(-G'(\beta_1))\mbox{   or  }4\pi/(G'(\beta_2)),\]
which are incompatible:
\[G'(\beta_1)+G'(\beta_2)=-\alpha(\beta_1-\beta_2)^2\neq 0.\]

\subsection{Pleba\'nski-Demia\'nski metrics}
These metrics, \cite{PD}, provide a useful step on the way to the Chen-Teo metric. In fact in \cite{CT} Chen and Teo describe their new metric as a bridge between Pleba\'nski-Demia\'nski and a 3-node Gibbons-Hawking metric. They can be Lorentzian or Riemannian and are still type-D (or Hermitian) but in the Riemannian case now ALE not ALF. A cosmological constant can be included, but we'll omit that.

Take the Riemannian version of the metric as given in \cite{CT}:
\be\label{pd1}
g=A^2dp^2+B^2dq^2+(\alpha d\tau+\beta d\phi)^2+(\lambda d\tau+\mu d\phi)^2,\ee
where
\[A^2=\frac{(1-p^2q^2)}{(p-q)^2P(p)},\;B^2=-\frac{(1-p^2q^2)}{(p-q)^2Q(q)}\]
\[\alpha^2=\frac{q^4P}{(p-q)^2(1-p^2q^2)},\;\beta^2=\frac{P}{(p-q)^2(1-p^2q^2)},\;\alpha\beta<0\]
\[\lambda^2=-\frac{Q}{(p-q)^2(1-p^2q^2)},\;\mu^2=-\frac{p^4Q}{(p-q)^2(1-p^2q^2)},\;\lambda\mu<0\]
and
\be\label{pd2}P(p)=\gamma+2np-\epsilon p^2+2mp^3+\gamma p^4,\;Q(q)=\gamma+2nq-\epsilon q^2+2mq^3+\gamma q^4.\ee
Although $P,Q$  are the same function, we'll adopt the convention that $P$ always means $P(p)$ and $Q$ always means $Q(q)$. There will be space-time singularities where $p^2q^2=1$, and points with $p=q$ must be at infinity. Very roughly, we can think of $m$ and $n$ as mass and nut-parameter, but $\epsilon$ is not charge (since the solution is vacuum).

For the signature to be Riemannian, we need $P$ positive and $Q$ negative (w.l.o.g. we can suppose $\gamma>0$). Thus, similarly to the $C$-metric, we need $p$ to lie between two roots of $P$ in an interval, say $I_1$, with $P$ positive, and $q$ to lie between two roots of $Q$ in an interval, say $I_2$, with $Q$ negative. We'll suppose that $P$ has four real roots, ordered $-\infty<p_1<p_2<0<p_3<p_4<\infty$ so there are three possibilities for $I_1$, namely $(-\infty,p_1),(p_2,p_3)$ or $(p_4,\infty)$, and two possibilities for $I_2$, namely $(p_1,p_2)$ or $(p_3,p_4)$, giving six possible rectangles altogether. The rectangles with a semi-infinite side can be ruled out as inevitably they will contain points with $p^2q^2=1$.  Therefore we want $p_2\leq p\leq p_3$ and either $p_1\leq q\leq p_2$ or $p_3\leq q\leq p_4$. Clearly points with $p=q$ are at infinity and one corner of either choice of the rectangle $I_1\times I_2$  lies on the line $p=q$, and is at infinity. The other corners give three nodes. The constants in $P$ must be chosen so that the singularities which occur at $p^2q^2=1$ don't meet the chosen rectangle, which we'll assume can be arranged for at least one of the two rectangles.

With the conventions of (\ref{1}) we have
\[U=\alpha^2+\lambda^2,\;V=\alpha\beta+\lambda\mu,\;W=\beta^2+\mu^2,\]
so that
\[r^2=UW-V^2=(\alpha\mu-\beta\lambda)^2\]
from which we may take
\be\label{pd3}r=\frac{(-PQ)^{1/2}}{(p-q)^2}.\ee
The axis is now at $PQ=0$ and the rods are at zeroes of $P$ or of $Q$ with nodes at common zeroes of both.
 
After a calculation we find the function conjugate to $r$, up to an additive constant, is
\be\label{pd4}z=(-\gamma-n(p+q)+\epsilon pq-\gamma p^2q^2-mpq(p+q))/(p-q)^2.\ee
On a rod or section of the axis where $P=0$ we may add $P$ to the numerator of $z$ to find the simplification
\[z=[(n-\epsilon  p+2mp^2+\gamma p^3)+(mp+\gamma p^2)q]/(p-q),\]
so that, since $p$ is constant on this rod, $z$ is a M\"obius fiunction of $q$. Furthermore along this rod
\[\frac{dz}{dq}=\frac{1}{2(p-q)^2}\frac{dP}{dp}.\]
 There is a similar conclusion on a rod with $Q=0$, when $z$ is a M\"obius function of $p$:
\[z=[(-n+\epsilon q-2mq^2-\gamma q^3)-(mq+\gamma q^2)p]/(p-q),\]
and
\[\frac{dz}{dp}=\frac{1}{2(p-q)^2}\frac{dQ}{dq}.\]

From (\ref{pd3}) we see that there are nodes at the common zeroes of $P$ and $Q$ so these will be at the corners of $I_1\times I_2$ and there are two possibilties:
\begin{itemize}
\item going clockwise round the rectangle from $(p_3,p_3)$, the point at infinity, the nodes are at
\[(p,q) =(p_2,p_3),\; (p_2,p_4),\;(p_3,p_4),\]
\item with corresponding $z$-values
\[z_1=-\frac{\epsilon}{2}+m(p_2+p_3)+\frac{\gamma}{2}(p_2+p_3)^2,\]
\[z_2=-\frac{\epsilon}{2}+m(p_2+p_4)+\frac{\gamma}{2}(p_2+p_4)^2,\]
\[z_3=-\frac{\epsilon}{2}+m(p_3+p_4)+\frac{\gamma}{2}(p_3+p_4)^2,\]
\item  by inspection of $dz/dq$ and $dz/dp$ along the rods, one sees that clockwise round the rectangle is the sense of increasing $z$  
so the top rod is $p=p_3$ and the bottom rod is $q=p_3$,
\end{itemize}
 or
\begin{itemize}
\item  going clockwise round the rectangle from $(p_2,p_2)$, the point at infinity, the nodes are at
\[(p,q) =(p_3,p_2),\;(p_3,p_1),\;(p_2,p_1),\]
\item   with corresponding $z$-values
\[z_1=-\frac{\epsilon}{2}+m(p_3+p_2)+\frac{\gamma}{2}(p_3+p_2)^2,\]
\[z_2=-\frac{\epsilon}{2}+m(p_3+p_1)+\frac{\gamma}{2}(p_3+p_1)^2,\]
\[z_3=-\frac{\epsilon}{2}+m(p_2+p_1)+\frac{\gamma}{2}(p_2+p_1)^2,\]
\item by inspection of $dz/dq$ and $dz/dp$ along the rods, again clockwise round the rectangle is the sense of increasing $z$ 
so the top rod is $p=p_2$, and the bottom rod is $q=p_2$.
\end{itemize}

\medskip

After some work, we obtain the twist potential for $K=\partial \tau$ by solving
\[d\psi=*(K\wedge dK)\]
 as
\be\label{pd5}\psi=\frac{F(p,q)}{(p-q)(1-p^2q^2)}\ee
with
\be\label{pd6} F=-2\gamma q-2nq^2+2mpq^3+2\gamma p^2q^3.\ee

Now we have enough to write down $P_+$ but the details are messy and will appear elsewhere. It's fairly straightforward to get the algebraic form as


\[P_+=\left(\begin{array}{cc}
         \frac{q_1}{C_1}& \frac{C_2}{C_1}\\
        \frac{C_2}{C_1}& \frac{Q_1}{C_1}\\
\end{array}\right),\]
where $q_1$ is quadratic in $z$, $C_1,C_2$ are both cubic and $Q_1$ is quartic. Recall these metrics are AE/ALE but still this is not in the form that one expects from (\ref{e8}). However one has the transformation in (\ref{gh5}) available and this can be used to remove the two highest terms in $C_2$, leaving it linear, say

\be\label{pd7}P_+=\left(\begin{array}{cc}
         \frac{q_1}{C_1}& \frac{\ell}{C_1}\\
        \frac{\ell}{C_1}& \frac{\widetilde{Q}_1}{C_1}\\
\end{array}\right),\ee
with linear $\ell(z)$, and $\widetilde{Q}_1$ a different quartic from the previous $Q_1$. Establishing this algebraic form for $P_+$ is the main aim of this subsection so we'll stop here and return to this form for $P_+$ in Section 5.3.

\subsubsection{Conicality}
There's no hope of removing all the conical singularities. The four rods have four different constant combinations of $\partial_\tau,\partial_\phi$ vanishing, in order from the bottom
\[p_3^2\p_\tau-\p_\phi,\;\p_\tau-p_4^2\p_\phi,\;p_2^2\p_\tau-\p_\phi,\;\p_\tau-p_3^2\p_\phi.\]

\subsection{The Chen-Teo metric}
Following \cite{CT}, the metric is
\be\label{ct1}
g=\frac{kH}{(x-y)^3}\left(\frac{dx^2}{X}-\frac{dy^2}{Y}- \frac{XY}{kF}d\phi^2\right)+\frac{(Fd\tau+Gd\phi)^2}{(x-y)HF},\ee
where
\bea\label{ct2}
X=a_0+a_1x+a_2x^2+a_3x^3+a_4x^4\\\label{ct3}
Y=a_0+a_1y+a_2y^2+a_3y^3+a_4y^4\\\label{ct4}
F=x^2Y-y^2X\\\label{ct5}
H=(\nu x+y)[(\nu x-y)(a_1-a_3xy)-2(1-\nu)(a_0-a_4x^2y^2)]\\\label{ct6}
G=(\nu^2a_0+2\nu a_3y^3+2\nu a_4y^4-a_4y^4)X+(a_0-2\nu a_0-2\nu a_1x-\nu^2a_4x^4)Y.
\eea
At this point there are seven parameters: $k$ which we mentioned in Section 2.1, $\nu$, and the five coefficients $a_0,\ldots,a_4$ of the polynomials $X$ and $Y$. We assume that $a_4>0$ and that $X$ (and therefore also $Y$) has four real roots, $x_i,i=1,\ldots,4$ with 
\[-\infty<x_1<x_2<0<x_3<x_4<\infty\]
and that $x$ is confined to an interval with $X$ positive, while $y$ is confined to an interval with $Y$ negative. As with the Pleba\'nski-Demia\'nski metric in the previous section, this leads to two possible rectangles in the $(x,y)$-plane, each with one vertex on the line $x=y$, which will be the point at infinity. In \cite{CT} we learn that we must take the rectangle $\{x_2\leq x\leq x_3,\;x_1\leq y\leq x_2\}$ to obtain a solution without curvature singularities (though conical singularities have yet to be eliminated). 

For the Weyl-Papapetrou coordinates, from \cite{CT}
\be\label{ct7}
\rho=\frac{(-XY)^{1/2}}{(x-y)^2},\;z=\frac{2(a_0+a_2xy+a_4x^2y^2)+(x+y)(a_1+a_3xy)}{2(x-y)^2},\ee
so that the axis is at $XY=0$, which is the boundary of whichever coordinate rectangle is chosen, and there are three nodes corresponding to the three corners of this rectangle not at infinity.

On $y=x_2$  claim
\[z=((a_1+2a_2x_2+2a_3x_2^2+2a_4x_2^3)+x(a_3x_2+2a_4x_2^2))(2(x-x_2))^{-1},\]
which is M\"obius in $x$ and positive, tending to infinity as $x$ tends from above to $x_2$, so this is the top rod.

On $x=x_2$ claim
\[z=((a_1+2a_2x_2+2a_3x_2^2+2a_4x_2^3)+y(a_3x_2+2a_4x_2^2))(2(y-x_2))^{-1},\]
which is M\"obius in $y$ and negative, tending to minus infinity as $y$ tends from below to $x_2$, so this is the bottom rod.

The topmost node is where the top rod arrives at $x_3$ so this is at
\[z=z_3= -\frac12(a_2+a_3(x_2+x_3)+a_4(x_2+x_3)^2) . \]
The bottom-most node is at $(x,y)=(x_2,x_1)$ so that
\[z=z_1= -\frac12(a_2+a_3(x_1+x_2)+a_4(x_1+x_2)^2) .  \]
(and it's easily checked that $z_3>z_1$).

For the Killing vector $K=\partial_\tau$ the twist potential $\psi$ solving
\[d\psi=*(K\wedge dK)\]
can be taken to be
\be\label{ct8}\psi=H^{-1}[a_0y+(2\nu-1)a_0x+(a_1\nu-a_4y^3)x^2+(-\nu a_3y-2\nu a_4y^2+a_4y^2)x^3].\ee
We've done enough to obtain $J'$ and $P_+$ but the details are complicated and can be found in \cite{DT24}. The basic shape, given ALF and three nodes, is
\be\label{ct9}
P_+=\Delta^{-1}\left(\begin{array}{cc}
         C_1& \ell\\
        \ell& C_2\\
\end{array}\right),\ee
with nodes at $z=a_1,a_2,a_3$, $\Delta(z)=(z-a_1)(z-a_2)(z-a_3)$, $C_i(z)$ cubic and $\ell(z)$ linear. Again we'll revisit this in Section 5.3.

\section{The inverse problem}
This is the problem of obtaining $P$ from the data, i.e. from the rod structure and the asymptotics,  rather than from the metric. It was solved for four-dimensional AF with 2 nodes in \cite{F,FW} and for some five-dimensional cases in \cite{NM2}.

Let us first note some general properties of $P(z)$ for a regular toric Ricci-flat metric with one of the asymptotic forms of AF/ALF/ALE: I claim it can't have poles off the axis as these would correspond to space-time singularities, and it can only have first-order poles on the axis at the nodes, by regularity. Justification of these claims will appear in a forthcoming article. If the nodes are located at $z=a_i,i=1,\ldots,n$ then every entry in $P$ is a rational function with denominator 
\[\Delta:=\Pi_1^n(z-a_i).\]
Furthermore the degrees of the numerators are fixed by the asymptotic conditions.

\medskip

I've made no mention of this, but it is known how to recognise (and therefore attempt to remove) conical singularities in this formalism: see p121 in \cite{F}. 

\medskip

If we can deduce $P(z)$ from the data then the metric can be obtained by splitting $P(z)$, and an algorithm for doinng this will appear in the near future.  Solving this inverse problem was part of the motivation of \cite{F,FW} as this was seen as a new approach to black hole uniqueness. The general technique, then, motivates the ansatz for AF/ALF
\[P_+=\frac{1}{\Delta}\left(\begin{array}{cc}
         P_1(z)& Q(z)\\
        Q(z)& P_2(z)\\
\end{array}\right),\]
where $P_1,P_2$ are polynomials of degree $n$ with leading coefficients $+1$ and $-1$ respectively and next coefficients $2M$ for both, and $Q(z)$ is a polynomial of degree $n-1$ with leading coefficient $2N$ and next coefficient $-2L$, while for ALE we take the ansatz
\[P_+=\frac{1}{\Delta}\left(\begin{array}{cc}
         P(z)& Q(z)\\
        Q(z)& S(z)\\
\end{array}\right),\]
where $P(z)$ is a polynomial of degree $n-1$ with leading coefficient $1/2$ and next coefficient $-M/2$, $Q(z)$ is a polynomial of degree $n-2$ with leading coefficient $L/2$, and $S(z)$ is a polynomial of degree $n+1$ with leading coefficient $-2$ and next coefficient $-2M$. Then the unit determinant condition will fix more of the free constants.

We'll order the discussion by the number of nodes. Zero nodes leads to flat space in the AF/ALF case and is inconsistent with AE/ALE so the simplest case to start with is a single node.
\subsection{One node}
For AF/ALF this case has only the SDTN solution from (\ref{tn8}):
\[P_+=
\frac{1}{z}\left(\begin{array}{cc}
         z+2m& -1\\
        -1& 0\\
\end{array}\right),\]
which we can transform according to (\ref{gh3}) with 
\[C=\left(\begin{array}{cc}
         1& 0\\
        1& 1\\
\end{array}\right),\]
to obtain 
\[P_+=\frac{1}{z}\left(\begin{array}{cc}
         z+2m& 2m\\
        2m& -z+2m\\
\end{array}\right),\]
and this takes the form of the ansatz but with $N=m$, and is the unique possibility.

For AE/ALE we only have $P_+$ as in (\ref{ee1}), and correspondingly the flat metric on $\mathbb{E}^4$.
\subsection{Two nodes}
We'll always place the nodes at $z=\pm \sigma$.

Then for AF/ALF the ansatz is
\[P_+=\frac{1}{z^2-\sigma^2}\left(\begin{array}{cc}
         z^2+2mz+A& 2Nz-2L\\
        2Nz-2L& -z^2+2mz+B\\
\end{array}\right).\]
The determinant condition branches into two cases:
\begin{itemize}
\item $m^2=N^2$: then the metric takes the Gibbons-Hawking form and the potential is for unequal masses at $z=\pm\sigma$. The solution isn't regular unless the masses are equal, when it's the two-centre SDTN.
\item $m^2\neq N^2$:  this is Kerr-Taub-bolt as in (\ref{kn2}); in particular it includes Kerr when $N=0$. This is the Riemannian counterpart of Black-hole Uniqueness: two nodes, Lorentzian and AF implies Kerr.
\end{itemize}

\medskip

With AE/ALE and two nodes, we can locate the nodes at $z=\pm\sigma$ as before and the asymptotics (\ref{e8}) force the ansatz:
\be\label{i1}
P_+=\frac{1}{(z-\sigma)(z+\sigma)}\left(\begin{array}{cc}
         \ell(z)& -2L\\
        -2L& c(z)\\
\end{array}\right),\ee
where $\ell(z)$ is linear, $L$ is constant, and $c(z)$ is cubic with some coefficients fixed by (\ref{e8}):
\[\ell(z)=\frac12(z-M),\;\;c(z)=-(2z^3+2Mz^2+2Az+2B),\]
where $M,L$ are the asymptotic quantities as in (\ref{e8}). The constants $A,B$ are to be fixed by the determinant condition. This gives
\[A=M^2-2\sigma^2,\;\;B=MA,\;\;\sigma^2=M^2\pm2L,\]
where we have solved a quadratic for $\sigma^2$ (which must be positive of course).

It turns out that this $P_+$ arises by appying the transformation in (\ref{gh5}) to the matric $P_+$ in (\ref{ba3}): take
\[\alpha=(2b)^{1/2},\;\beta=(2/b)^{1/2},\;\gamma=a\sigma(2/b^3)^{1/2},\]
then $P_+$ from (\ref{ba3}) precisely becomes $P_+$ of (\ref{i1}) with
\[M=a\sigma/b,\;L=\frac{\sigma^2(a^2-b^2)}{2b^2}.\]
Therefore the general ALE two-node solution is the Baza\u{i}kin solution (\ref{ba3}), which includes the Eguchi-Hanson metric.

\subsection{Three nodes}

Now for three nodes, suppose they are at $z=\alpha,\beta$ and $\gamma$, where w.l.o.g. $\alpha+\beta+\gamma=0$, and set 
\[\Delta=(z-\alpha)(z-\beta)(z-\gamma).\]
The ansatz for an AF/ALF solution is
\[P_+(z)=\Delta^{-1}\left(\begin{array}{cc}
         C_1(z)& \ell(z)\\
        \ell(z)& C_2(z)\\
\end{array}\right),\]
with
\[C_1(z)=z^3+2Mz^2+Az+B,\]
\[\ell(z)=2L(z-C)\]
\[C_2(2)=-z^3+2Mz^2+Dz+E,\]
with the determinant condition
\[C_1(z)C_2(z)-(q(z))^2=-\Delta^2.\]
This is analysed in \cite{DT24} and has been the subject of further work, to appear \cite{MD}. In brief, the solution space for $P_+$ branches: one branch is the 3-node multi-SDTN solutions and the other is the 5-parameter Chen-Teo solution of \cite{CT}.

\medskip

The ansatz for an ALE solution is
\[P_+(z)=\Delta^{-1}\left(\begin{array}{cc}
         q(z)& \ell(z)\\
        \ell(z)& Q(z)\\
\end{array}\right),\]
with
\[q(z)=\frac12(z^2+2Mz+2A),\]
\[\ell(z)=2L(z-B)\]
\[Q(2)=-2(z^4-2Mz^3+Cz^2+Dz+E),\]
and $A,B,C,D,E$ free. The determinant condition is
\[q(z)Q(z)-(\ell(z))^2=-\Delta^2.\]
It can be solved for all quantities in terms of $A$, leaving a quartic for $A$. Again the solution space branches: one branch is the 3-node Eguchi-Hanson solutions and the other is the Pleba\'nski-Demia\'nski solutions.

%
%
%
\section*{Appendix 1: Double poles}
In a forth-coming paper it will be shown that, for a complete Riemannian metric and AF/ALF/ALE asymptotics, $P_+$ can only have a finite number of simple poles, located at the nodes. The story is different in Lorentzian signature where a double pole can be the sign of a regular but extremal horizon. One can see this by taking $\sigma$ to zero in (\ref{k10}), the $P_+$ for Lorentzian Kerr, when the result has a double pole at the remaining node but is otherwisw unremarkable. If one takes the same limit, i.e. of vanishing $\sigma$, in (\ref{k7}) for Riemannian Kerr then $m$ and $a$ both go to zero and the result is (one choice of) the $P_+$ for flat space.

If we try this for Taub-NUT, for the Lorentzian case (\ref{tn11}) as $\sigma$ goes to zero, both parameters $m$ and $\ell$ go to zero, leading to flat Minkowski space, while for the Riemannian case $N\rightarrow m$ and the limiting $P_+$ is
\[\left(\begin{array}{cc}
    1+2m/z & 2m/z \\
    2m/z & -1+2m/z\\
\end{array}\right),\]
which is actually of Gibbons-Hawking form as can be seen by the transformation
\[P_+\rightarrow \left(\begin{array}{cc}
    1 & 0\\
    -1 & 1\\
\end{array}\right)\left(\begin{array}{cc}
    1+2m/z & 2m/z \\
    2m/z & -1+2m/z\\
\end{array}\right)\left(\begin{array}{cc}
    1& -1 \\
    0& 1\\
\end{array}\right)=\left(\begin{array}{cc}
    1+2m/z & -1\\
    -1 & 0\\
\end{array}\right),\]
so this is (\ref{tn8}), the Self-dual Taub-NUT in Gibbons-Hawking form.

\medskip

One can readily construct Gibbons-Hawking examples with a double pole on the axis e.g. take for $V$ a simple dipole:
\[V=\frac{z}{(r^2+z^2)^{3/2}}\]
so that $V(0,z)=z^{-2}$ and $P$ in (\ref{gh3}) has a double pole at the origin. However the 4-metric is singular there: as we know anything but a simple pole in $V$ gives a space-time singularity. This follows from the splitting process, or from the observation that the only regular GH metrics are multi-Taub-NUT or multi-Eguchi-Hanson -- this dipole can arise as a limit of equal but oppositely charged monopoles, but that is not a regular solution.

\medskip

We haven't considered the stationary, axisymmetric Tomimatsu-Sato solutions, \cite{TS}, since they have naked singularities. If we consider the Tomimatsu-Sato solution with the parameter $\delta=2$ and Lorentzian signature, it depends on parameters $p,q$ with $p^2+q^2=1$. We can take the details from \cite{JI} where the twist potential $\psi$ is also calculated. On the top axis, $z\geq 1$, in the Lorentzian form we find
\[ P_+= \left(\begin{array}{cc}
    \frac{Q(z)}{p^2(z^2-1)^2} & \frac{8qz^2}{p^2(z^2-1)^2} \\
    \frac{8qz^2}{p^2(z^2-1)^2} & \frac{p^4(z^2-1)^4+64q^2z^4}{p^2(z^2-1)^2Q(z)}\\
\end{array}\right),\]
with
\[Q(z)=p^2(z^2-1)^2+4pz(z^2+1)+8z^2=(p(z^2+1)+2z)^2+4q^2z^2.\]
We've written $Q$ as a sum of squares, to show that it doesn't have real zeroes in $z$, but it's obtained from the norm of the $\partial_t$-Killing vector which is singular on a ring in the equatorial plane, giving a naked singularity, \cite{GRC}.

The $P_+$ matrix, and indeed the metric, has AF asymptotics with $m=2/p, L=-4q/p^2$, but $P_+$ has double poles at $z=\pm 1$ and this metric does not have a smooth extremal horizon: these double poles are not the sign of a smooth space-time feature.

\section*{Appendix 2: Splitting}
This is an important part of the Ward construction \cite{W2} and it has a particular form in this special case with two symmetries -- see \cite{W}, \cite{FW}. These two references differ in their details so I'm going to follow the latter and give just an abbreviated account of it. I'll start with a motivating example, which is classical\footnote{In the sense of old enough not to need a reference.}: suppose $F(z)$ is holomorphic in some part of the complex $z$-plane including at least part of the real axis. Substitute $z-\frac12r(\zeta-\zeta^{-1})$ for $z$ and expand $F$ as a Laurent series in $\zeta$:
\[F(z-\frac12r(\zeta-\zeta^{-1}))=\Sigma_{-\infty}^{\infty}a_n(r,z)\zeta^n,\]
then differentiating the series term by term we find that
\[\Sigma_{-\infty}^\infty[ra_{n,rr}+a_{n,r}+ra_{n,zz}-\frac{n^2}{r}a_n]\zeta^n=0.\]
In particular therefore $a_0(r,z)$ is harmonic and this term can be extracted by a contour integral:
\[F(r,z):=a_0(r,z)=\frac{1}{2\pi i}\oint F(z-\frac12r(\zeta-\zeta^{-1}))\frac{d\zeta}{\zeta}.\]
Furthermore, evaluating at $r=0$ we find $F(0,z)=F(z)$ so this is the harmonic function with value $F(z)$ on the axis.

We connect to the patching matrix following \cite{FW}: with patching matrix $P(z)$ substitute for $z$ as before and split by factorising:
\[P(z-\frac12r(\zeta-\zeta^{-1}))=H(r,z;\zeta)\widehat{H}^{-1}(r,z;\zeta),\]
where $\widehat{H}(r,z;\zeta)$ is holomorphic in $\zeta$ for $|\zeta|\geq 1$, including $\infty$, and $H(r,z;\zeta)$ is holomorphic in $\zeta$ for $|\zeta|\leq 1$ (this is always possible). Now set
\be\label{sp1}J'(r,z)=H(r,z;0)\widehat{H}^{-1}(r,z;\infty),\ee
then this is the desired Ernst potential from $P(z)$ and satisfies the Yang equation (\ref{11}).

\medskip

For a simple example, consider (\ref{16}) with a static metric, so that $\psi=0$ and by (\ref{p1})
\[P(z)= \left(\begin{array}{cc}
   (f(0,z))^{-1} & 0 \\
    0 & f(0,z)\\
\end{array}\right),\]
and replace $z$ by $z-\frac12r(\zeta-\zeta^{-1})$. Evidently $H,\widehat{H}$ can be taken to be diagonal:
\[\widehat{H}=\left(\begin{array}{cc}
   e^{\hat{h}} & 0 \\
    0 & e^{-\hat{h}}\\
\end{array}\right),\;H=\left(\begin{array}{cc}
   e^{-h} & 0 \\
    0 & e^h\\
\end{array}\right)\]
and so we need
\[\hat{h}(r,z;\zeta)+h(r,z;\zeta)=\log f(z-\frac12r(\zeta-\zeta^{-1})),\]
which we solve by expanding the r.h.s. as a Laurent series in $\zeta$, as before. Call the $O(0)$ term $a_0(r,z)$ then the prescription (\ref{sp1}) gives
\[J'=\left(\begin{array}{cc}
   e^{-a_0(r,z)} & 0 \\
    0 & e^{a_0(r,z)}\\
\end{array}\right).\]
We recover $J(r,z)$ from (\ref{14}) and the metric by solving (\ref{9},{10}) with $\omega=0$. These are the Weyl solutions -- see Section 4.6 or \cite{ES}.

\medskip

We can use the same trick to split the patching matrix for the Gibbons-Hawking solutions (\ref{gh3}): put $z-\frac12r(\zeta-\zeta^{-1})$ for $z$ and consider
\[P(z-\frac12r(\zeta-\zeta^{-1}))=\left(\begin{array}{cc}
   F(z-\frac12r(\zeta-\zeta^{-1})) & -1 \\
    -1 & 0\\
\end{array}\right) =H(r,z;\zeta)\widehat{H}^{-1}(r,z;\zeta),\]
where $F(z)=V(0,z)$.
A shrewd choice is
\[  H(r,z;\zeta)=\left(\begin{array}{cc}
   1 & h \\
    0 &-1\\
\end{array}\right),\;\;\widehat{H}(r,z;\zeta)  =\left(\begin{array}{cc}
   0 & 1 \\
    -1 & \hat{h}\\
\end{array}\right),\]
and then we just need
\[ F(z-\frac12r(\zeta-\zeta^{-1}))=\hat{h}+h\]
so that, with the Laurent series $F=\Sigma_{-\infty}^{\infty}a_n(r,z)\zeta^n$, we may take
\[\hat{h}=\Sigma_{-\infty}^{-1}a_n(r,z)\zeta^n,\;\;h=\Sigma_0^{\infty}a_n(r,z)\zeta^n,\]
and obtain
\[J'(r,z)=H(r,z;0)\widehat{H}^{-1}(r,z;\infty)=\left(\begin{array}{cc}
   a_0(r,z) & -1 \\
    -1 & 0\\
\end{array}\right),\]
which is correct since w.l.o.g. $\psi=f=V^{-1}$ in (\ref{15}), so that $V(r,z)=a_0(r,z)$ is harmonic.

\medskip

A general splitting algorithm is under development and should be forthcoming..

\end{document}